\documentclass{amsart}

\usepackage{amssymb}
\usepackage[all]{xy}
\usepackage{hyperref}

\newtheorem{thm}{Theorem}

\newtheorem{lem}[thm]{Lemma}
\newtheorem{cor}[thm]{Corollary}

\newtheorem{prop}[thm]{Proposition}

\newtheorem{conj}[thm]{Conjecture}
   
\theoremstyle{definition}
\newtheorem{defn}[thm]{Definition}

\newtheorem{say}[thm]{}
\newtheorem{exmp}[thm]{Example}

\newtheorem{rem}[thm]{Remark}          

\newtheorem*{ack}{Acknowledgments}      

\newtheorem{defn-thm}[thm]{Definition--Theorem}  
\newtheorem{defn-lem}[thm]{Definition--Lemma}  

\theoremstyle{remark}

\renewcommand{\c}[0]{{\mathbb C}}  

\renewcommand{\o}[0]{{\mathcal O}} 
\newcommand{\z}[0]{{\mathbb Z}}
\newcommand{\n}[0]{{\mathbb N}}
\renewcommand{\r}[0]{{\mathbb R}} 

\renewcommand{\a}[0]{{\mathbb A}}

\newcommand{\p}[0]{{\mathbb P}}
\newcommand{\f}[0]{{\mathbb F}}
\newcommand{\q}[0]{{\mathbb Q}}
\newcommand{\map}[0]{\dasharrow}
\newcommand{\qtq}[1]{\quad\mbox{#1}\quad}
\newcommand{\spec}[0]{\operatorname{Spec}}

\newcommand{\gal}[0]{\operatorname{Gal}}

\newcommand{\mult}[0]{\operatorname{mult}}

\newcommand{\proj}[0]{\operatorname{Proj}}

\newcommand{\coker}[0]{\operatorname{coker}}

\newcommand{\aut}[0]{\operatorname{Aut}}

\newcommand{\sing}[0]{\operatorname{Sing}}

\newcommand{\chr}[0]{\operatorname{char}}

\newcommand{\gr}[0]{\operatorname{gr}}

\newcommand{\rdown}[1]{\lfloor{#1}\rfloor}

\newcommand{\diag}[0]{\operatorname{diag}}

\newcommand{\tsum}[0]{\textstyle{\sum}}

\newcommand{\bir}[0]{\operatorname{Bir}}

\def\into{\DOTSB\lhook\joinrel\to}

\def\loccoh#1.#2.#3.#4.{H^{#1}_{#2}(#3,#4)}

\DeclareMathAlphabet{\mathchanc}{OT1}{pzc}%
                                {m}{it}

\newcommand{\gm}[0]{{\mathbb G}_m}
\newcommand{\ga}[0]{{\mathbb G}_a}

\newcommand{\GL}{\mathrm{GL}}
\newcommand{\PGL}{\mathrm{PGL}}
\newcommand{\PSL}{\mathrm{PSL}}
\newcommand{\SU}{\mathrm{SU}}
\newcommand{\PSU}{\mathrm{PSU}}
\newcommand{\SL}{\mathrm{SL}}

\newcommand{\SO}{\mathrm{SO}}
\newcommand{\PSO}{\mathrm{PSO}}
\newcommand{\PO}{\mathrm{PO}}
\newcommand{\OO}{\mathrm{O}}

\newcommand{\GO}{\mathrm{GO}}
\newcommand{\PSp}{\mathrm{PSp}}

\newcommand{\sym}[0]{\operatorname{Sym}}

\usepackage[all]{xy}\xyoption{dvips}

\newcommand{\tprod}[0]{\textstyle{\prod}}

\newcommand{\auts}[0]{\operatorname{\mathbf{Aut}}}
\newcommand{\autu}[0]{\operatorname{Aut}^{\rm u}} 
  
\begin{document}
\bibliographystyle{amsalpha}


\title[Rings of invariants]{Automorphisms and twisted forms \\ of rings of invariants}
        \author{J\'anos Koll\'ar}

        \begin{abstract}  Let  $G\subset \GL_n(k)$ be a  finite  subgroup
          and $k[x_1,\dots, x_n]^G\subset k[x_1,\dots, x_n]$ its ring of invariants. We show that, in many cases,  the automorphism group of $k[x_1,\dots, x_n]^G$ is $k^\times$.
        \end{abstract}

\maketitle

Let  $k$ be a field and $G\subset \GL_n(k)$  a  finite  subgroup.
Then $G$ also acts on the polynomial ring $k[x_1,\dots, x_n]$ and on the power series  ring $k[[x_1,\dots, x_n]]$, which we
abbreviate as $k[{\mathbf x}]$ (resp.\  $k[[{\mathbf x}]]$). The ring of invariants is denoted by $k[{\mathbf x}]^G$  (resp.\  $k[[{\mathbf x}]]^G$). Our aim is to describe the automorphism groups of these rings.

Let  $N(G):=N(G,\GL_n(k))$ denote the normalizer of $G$ in $\GL_n(k)$.
Then the $N(G)$-action on $k[{\mathbf x}]$ descends to an
$N(G)/G$-action on $k[{\mathbf x}]^G$.
We prove that these automorphisms of $k[{\mathbf x}]^G$  and  $k[[{\mathbf x}]]^G$ give the largest reductive quotient/subgroup of the automorphism group.

The normalizer always contains the scalars  $k^\times$. They act
on  homogeneous polynomials as $p({\mathbf x})\mapsto \lambda^{\deg p}p({\mathbf x})$, but  we get a sometimes larger  $k^\times$-action given as
$p({\mathbf x})\mapsto \lambda^{(\deg p)/e}p({\mathbf x})$, where
$e$ is the $\gcd$ of the degrees of all homogeneous $G$-invariants.


Other automorphisms are not always easy to write down, and it turns out that, depending on $G$,   $\aut \bigl(k[{\mathbf x}]^G\bigr)$ is sometimes infinite dimensional, and sometimes 1 dimensional.
We give criteria for both extremes, but  a complete answer remains conjectural.

More generally, we consider rings of invariants by  finite, geometrically reduced subgroup schemes   $G\subset \GL_n$  over $k$.


If $\chr k$ does not divide $|G|$, then,   by the Chevalley-Shephard-Todd theorem \cite{chevalley}, $k[{\mathbf x}]^G\cong k[{\mathbf y}]^H$ for a suitable $H$-action, where  $H$ is the quotient of $G$ by the subgroup generated by all 
 pseudo-reflections, that is, elements that fix a hyperplane pointwise.
Thus we can harmlessly assume that  $G$ does not contain any 
pseudo-reflections. The situation is less clear if 
$\chr k$  divides $|G|$.

Our first result computes $\aut \bigl(k[{\mathbf x}]^G\bigr)$ in terms of
$G\subset \aut\bigl(k[{\mathbf x}]\bigr)$.

\begin{thm} \label{main.k.thm.i}
  Let $k$ be a field and $G\subset \GL_n$ a  finite, geometrically reduced subgroup scheme  without pseudo-reflections. Assume that
  \begin{enumerate}
  \item either $R=k[[{\mathbf x}]]$,
  \item or  $R=k[{\mathbf x}]$, $\chr k\nmid |G|$, and  $0\in k^n$ is the only $G$-invariant vector.
    \end{enumerate}
  Then there is a split exact sequence 
  $$
  1\to C\bigl(G, \autu(R)\bigr) \to  \aut (R^G)\to
  \bigl({\bf N}(G,\GL_n) /G\bigr)(k)\to 1.
  \eqno{(\ref{main.k.thm.i}.3)}
  $$
\end{thm}

Here $\autu(R)\subset \aut(R)$ denotes the subgroup of those automorphisms that stabilize  the maximal ideal $m=(x_1,\dots, x_n)$ and act trivially on $m/m^2$.  It is the group of $k$-points of an infinite dimensional ind-group
\cite{MR485898, MR607583, https://doi.org/10.48550/arxiv.1809.04175}, but the 
ind-group structure plays no role for us.

$C(G, *)$ denotes the centralizer, the   subgroup of those elements of $*$  that commute with $G$.
If the origin is the only $G$-fixed point, then  the maximal ideal $m$  is also fixed by the normalizer
$N\bigl(G, \aut(R)\bigr)$, so
$N\bigl(G, \aut(R)\bigr)=N\bigl(G, \aut(R,m)\bigr)$.
Since $\autu(R) $ is a normal subgroup of $\aut(R,m)$, and its  intersection with $\GL_n$ is the identity, the normalizer 
$N\bigl(G, \autu(R)\bigr)$  equals the centralizer
$C\bigl(G, \autu(R)\bigr)$.

On the right in (\ref{main.k.thm.i}.3),  ${\bf N}(G,\GL_n)$ denotes the normalizer of $G$ in $\GL_n$, as an algebraic group over $k$. We take the quotient by $G$, as an algebraic group, and then
the  $k$-points of the quotient. 
 It is usually  larger than
 $ N(G,\GL_n(k)) /G(k)$, but the two agree if either $k$ is separably closed, or if $G$ acts irreducibly, and its center is trivial;  see  Examples~\ref{An.long.exmp}--\ref{C=bfC.exmp}.
 In any case, this part is  quite well understood.

 On the left in (\ref{main.k.thm.i}.3)
$C\bigl(G, \autu(k[[{\mathbf x}]])\bigr)$ is always infinite dimensional, see Example~\ref{uip.inf.ps.exmp}.
 By contrast, $C\bigl(G, \autu(k[{\mathbf x}])\bigr)$
is  sometimes trivial.
 It is    infinite dimensional,  whenever $G$ is diagonalizable over $k$, see Example~\ref{reducible.exmp.0}. The same holds for most reducible representations by Example~\ref{reducible.exmp.1}.
However, for irreducible representations, the following could be true.

\begin{conj} \label{centr.conj} Let $k$ be a field and $G\subset \GL_n$ a  finite, geometrically reduced subgroup scheme that acts irreducibly over $k$. Then
  $ C\bigl(G, \autu(k[{\mathbf x}])\bigr)=\{1\}$.
  \end{conj}

We show this in dimension 2 for any $k$ in Theorem~\ref{2d.irred.thm.1},
and in dimension 3 for  $k=\c$ in Proposition~\ref{3d.irred.thm.1}. 
In higher dimensions,  we develop a notion of
 {\it birationally irreducible} representations (see Definition~\ref{bir.red.defn}) and  give a series of examples in (\ref{k-sz-fix.lem.c1}--\ref{pgl.q.rigid.prop}).
 (Very roughly speaking, $G$ is   birationally irreducible if the induced action
 on $\p^{n-1}_k$ is not birational to an action on a product.)
 We prove the following special case of Conjecture~\ref{centr.conj}.


\begin{thm}\label{bir.ired.act.cor.i}
  Let $G\subset \GL_n$ be a geometrically reduced and birationally irreducible  subgroup scheme.
  Then  the centralizer $C\bigl(G, \autu(k[{\mathbf x}])\bigr)$ is trivial.
\end{thm}

Combining the two theorems  gives the following.

\begin{cor}  \label{nd.irred.thm.i}
  Let $k$ be a field and $G\subset \GL_n$ a finite, birationally irreducible  subgroup  scheme without pseudo-reflections.  Assume that     $\chr k\nmid |G|$. Then
  $$
  \aut \bigl(k[{\mathbf x}]^G\bigr)=
  \bigl({\bf N}(G,\GL_n) /G\bigr)(k).  \hfill\qed
  $$
\end{cor}

The sequence (\ref{main.k.thm.i}.3) also allows us to determine the
 twisted forms of quotient singularities, where
 a {\it  twisted form} of a $k$-algebra $A$ is a $k$-algebra $A'$, such that
 ${\bar k}\otimes_kA\cong {\bar k}\otimes_k A'$.
 For illustrations see Examples~\ref{exmp.2}--\ref{exmp.1}.


 \begin{cor}\label{main.thm}
   Let $k$ be a perfect field with  algebraic closure ${\bar k}\supset k$ and
      $G\subset \GL_n$  a finite, geometrically reduced $k$-subgroup scheme without pseudo-reflections.
  Let $N\bigl(G,\GL_n({\bar k})\bigr)  \subset \GL_n({\bar k})$ be the normalizer of $G$.  Then there are natural bijections between the sets
  \begin{enumerate}
  \item twisted forms of $k[[{\mathbf x}]]^G$ (up to isomorphisms),
    \item  grading preserving twisted forms of $k[{\mathbf x}]^G$ (up to isomorphisms), and
\item $H^1\bigl(\gal({\bar k}/k), N\bigl(G,\GL_n({\bar k})\bigr)/G({\bar k})\bigr)$.
\end{enumerate}
  \end{cor}

 The proof of Theorem~\ref{main.k.thm.i} builds on \cite{MR3089030}. A key case is computed in \cite{https://doi.org/10.48550/arxiv.2210.12781}. The 
proof of Theorem~\ref{bir.ired.act.cor.i} relies on the method of
\cite{MR1288046}, as refined by \cite{MR2340971}.
See \cite{https://doi.org/10.48550/arxiv.2301.13040} for some related examples.
Corollary~\ref{main.thm} was inspired by \cite{bre-vis-2}.

 \begin{ack}  I thank  G.~Bresciani, Y.~Tschinkel, U.~Umirbaev, A.~Vistoli  and C.~Xu for   useful comments, and  P.~Tiep for giving the examples in   (\ref{no.hom.exmps}.4).
Partial  financial support    was provided  by  the NSF under grant number
DMS-1901855.
\end{ack}

\section{Examples}

We give a series of examples illustrating various possibilities for 
$\aut \bigl(k[{\mathbf x}]^G\bigr)$.

\begin{exmp} \label{An.long.exmp}
The automorphism groups of the $A_{n-1}$ singularities are quite interesting.
  Over a field $k$, write it as $(x^2-dy^2=z^n)\subset \a^3_k$ and assume that $\chr k\nmid 2n$.

  If $d$ is a square, it is equivalent to $(xy=z^n)$.
  This is the invariant ring of the $\mu_n$-action given by
  $(u,v)\mapsto (\epsilon u, \epsilon^{-1}v)$, where $\mu_n\subset \gm$ is the group of $n$th roots of unity.
  For any polynomial $p(x)$ we get an automorphism
  $$
  (x,y,z)\mapsto\Bigl(x, y+\frac{(z+xp(x))^n-z^n}{x}, z+xp(x)\Bigr),
  \eqno{(\ref{An.long.exmp}.1)}
    $$
  hence $\aut(xy=z^n)$ is infinite dimensional.

  If $d$ is not  a square, then $(x^2-dy^2=z^n)$ is the invariant ring $k[u,v]^G$, where $G$ is  the
  $n$-torsion subgroup of  $\SO (u^2-dv^2)$. This is not diagonalizable over $k$.    The connected component of its normalizer is
    $$
  \GO(u^2-dv^2)=
  \left\{
\begin{pmatrix}
  c & s\\
  ds & c
\end{pmatrix}
\colon c^2-ds^2\neq 0
\right\}\subset \GL_2.
\eqno{(\ref{An.long.exmp}.2)}
$$
Adding the reflection $(u,v)\mapsto (-u,v)$ gives the whole normalizer.

Using Corollary~\ref{nd.irred.thm.i} and Theorem~\ref{2d.irred.thm.1}, we  need to compute $\GO(u^2-dv^2)/G$. We have  2 cases. 
\medskip

{\it Case \ref{An.long.exmp}.3.} If $n=2r+1$, then
the reflection descends to  $(x,y)\mapsto (-x,y)$, and 
$$
\aut^\circ(x^2-dy^2=z^n)\cong  \GO(u^2-dv^2)/G\cong \GO(x^2-dy^2),
$$
where the quotient map is given by
$N\mapsto \det^{-r}N\cdot N^{2r+1}$,  with
$M\in \GO(x^2-dy^2)$ acting as
$$
\begin{pmatrix}
  x\\
  y\\
  z
\end{pmatrix}
\mapsto
\begin{pmatrix}
{\det}^rM\cdot M\cdot \begin{pmatrix}
  x\\
  y
\end{pmatrix}\\
\det M\cdot z
\end{pmatrix}.
$$

{\it Case \ref{An.long.exmp}.4.} If $n=2r$, then
the reflection descends to  $(x,y)\mapsto (x,-y)$, and 
$$
\aut^\circ(x^2-dy^2=z^{n+1})\cong \GO(u^2-dv^2)/G\cong \gm\times \OO(x^2-dy^2),
$$ where  the quotient map is given by
$N\mapsto \bigl(\det N, \det^{-r}N\cdot N^{2r}\bigr)$, with  $(\lambda, M)\in \gm\times \OO(x^2-dy^2)$ acting as
$$
\begin{pmatrix}
  x\\
  y\\
  z
\end{pmatrix}
\mapsto
\begin{pmatrix}
{\lambda}^r\cdot M\cdot \begin{pmatrix}
  x\\
  y
\end{pmatrix}\\
\lambda\cdot z
\end{pmatrix}.
$$

In particular, the ind-group  $\aut_{\c}(x^2+y^2=z^n)$ is
infinite dimensional, but its real form 
$\aut_{\r}(x^2+y^2=z^n)$ is 2-dimensional.
Similar examples are in \cite{MR2340971}.
(Note that for a finite dimensional group, the dimension is the same over $\r$ and $\c$.)
\end{exmp}

\begin{exmp}  \label{E8.long.exmp}
  Let $2\rm{I}\subset \SL_2(k)$ be the binary icosahedral group.
 Assume that  $\chr k\nmid |2\rm{I}|=120$.
  The ring of invariants is
   $k[x,y,z]/(x^2+y^3+z^5)$, where
  the degrees of $x,y,z$ are  $30, 20,12$. 

  Since $\rm{I}\subset \PGL_2$ is a maximal finite subgroup, its normalizer is itself. So ${\bf N}(2\rm{I},\GL_2)$ is  the product of $2\rm{I}$ and of 
  ${\bf C}(2\rm{I},\GL_2)$; the latter is the group of scalars  $\gm$.
  The diagonal matrix 
$\diag(c,c)$ acts on the invariants as
$$
\tau_c: (x,y,z)\mapsto  \bigl(c^{30}x, c^{20}y, c^{12}z\bigr).
$$
  However   $\diag(-1,-1)\in 2\rm{I}$, so
  the automorphism group is the $k$-points of $\gm/\{\pm 1\}$.
  As in  \cite{MR3089030} we get that 
  $$
  \aut\bigl(k[x,y,z]/(x^2+y^3+z^5)\bigr)=\bigl\{\tau'_c:
(x,y,z)\mapsto  (c^{15}x, c^{10}y, c^{6}z)\bigr\}\cong k^{\times}.
$$
\end{exmp}

\begin{exmp}\label{C=bfC.exmp}
More generally, 
  let $G\subset \GL_n$ be a   finite,  geometrically irreducible  subgroup 
  whose image  in $\PGL_n$ is maximal among finite subgroups.
  Then 
  ${\bf C}(G,\GL_n)\cong \gm$  surjects onto ${\bf N}(G,\GL_n)/G$, and so
  ${\bf N}(G,\GL_n)/G\cong \gm/(\mbox{center of $G$})\cong \gm$.

 If the assumptions of Corollary~\ref{nd.irred.thm.i} also hold, then
 $\aut\bigl(k[{\mathbf x}]^G\bigr)\cong k^\times$.  
   \end{exmp}

\begin{exmp}\label{reducible.exmp.0} Let $A$ be a commutative group scheme with a faithful representation on $k^n$ given by characters  $\chi_i:A\to \gm$.
  If $\chi_1, \dots, \chi_{n-1}$ generate the whole character group of $A$, then there are infinitely many $(m_1, \dots,m_n)\in\n^n$ such that $\chi_1^{m_1}\cdots\chi_{n-1}^{m_{n-1}}=\chi_n$.
  Then
  $$
  (x_1,\dots, x_n)\mapsto (x_1,\dots, x_{n-1}, x_n+cx_1^{m_1}\cdots x_{n-1}^{m_{n-1}})
  $$
  commutes with $A$, and descends to an automorphism of $k[{\mathbf x}]^A$. Thus
  $\aut\bigl(k[{\mathbf x}]^A\bigr)$ is infinite dimensional.

  If $\chi_1, \dots, \chi_{n-1}$ do not generate the whole character group, then
  there is an $a\in A$ such that $\bigl(\chi_1(a), \dots, \chi_{n}(a)\bigr)$
  is a pseudo-reflection.
  \end{exmp}

\begin{exmp}\label{reducible.exmp.1}
  Consider 2 representations  $G\to \GL(U)$ and $G\to \GL(V)$.
  Assume that there is a $G$-equivariant linear map $\phi: \sym^m(U)\to V$.  Let $p$ be a polynomial of degree $m$ on $U$. Then
  $
  (u,v)\mapsto \bigl(u, v+\phi(p(u))\bigr)
  $
  is an automorphism of $U+V$ that commutes with the diagonal $G$-action.

  If $G\into \GL(U)$,  then   there are nonzero  $\phi: \sym^{m}(U)\to V$ by
\cite{ko-ti}.
\end{exmp}

\begin{exmp}\label{uip.inf.ps.exmp}
  Consider a finite subgroup  $G\subset \GL_n$.
  Assume that there is a $G$-equivariant linear map $\phi: \sym^m(V)\to V$ for some $m\geq 2$.  Let $p$ be a polynomial of degree $m$ on $V$. Then
  $  v\mapsto v+\phi(p(v)) $
  is an automorphism of $k[[{\mathbf x}]]$ that commutes with the  $G$-action. (It is usually not  an automorphism of $k[{\mathbf x}]$.)

 For any $V$ there are infinitely many  $\phi: \sym^{m}(V)\to V$ by
 \cite{ko-ti}.
  \end{exmp}

\begin{exmp} Let $k$ be a field of characteristic 2 and consider the $\z/2$-action
  $$
  (x_1, y_1, x_2, y_2)\mapsto (x_1, y_1+x_1, x_2, y_2+x_2).
  $$
  The invariants are
  $$
  u_i:=x_i,\quad v_i:=y_i(y_i+x_i),\qtq{and} w:=x_1y_2+x_2y_1.
  $$
  They satisfy the equation
  $$
  w^2+u_1u_2w+u_1^2v_2+u_2^2v_1=0.
  $$
  Consider the automorphism $\phi_m$ that fixes $u_1, u_2, w$ and sends
$$
  (v_1, v_2)\mapsto (v_1+u_1^m, v_2+u_2^2u_1^{m-2}).
  $$
  Thus
  $\phi_m\bigl(y_1(y_1+x_1)\bigr)=y_1^2+y_1x_1+x_1^m$, and 
  the right had side is irreducible for $m\geq 3$. So there is no automorphism of $k[x_1, y_2, x_2, y_2]$ that sends $y_1(y_1+x_1)$ to $ y_1^2+y_1x_1+x_1^m$.

  Thus the assumption  $\chr k\nmid |G|$ is necessary in
  Theorems~\ref{main.k.thm.i}~and~\ref{main.u.k.thm}.
  It is not clear what happens for irreducible or birationally irreducible actions.
\end{exmp}

Commuting with pseudo-reflections is a strong restriction.

\begin{exmp}[Actions with  many pseudo-reflections]\label{man.ps.r.say}
  Any automorphism $\phi\in \aut(\a^n)$ that commutes with a pseudo-reflection fixes the invariant hyperplane. If there are $n$ such linearly independent hyperplanes  $(x_i=0)$, then  $\phi$ is an automorphism of
  $\gm^n\cong \a^n\setminus(x_1\cdots x_n=0)$. Note that the units on 
  $\gm^n$ are the monomials in the $x_i$, so $\aut(\gm^n)$ maps monomials to monomials.
  The only ones that  extend to  $\aut(\a^n)$  are
  $(x_1,\dots, x_n)\mapsto (c_1x_1,\dots, c_nx_n)$, where $c_i\in k^\times$.

  If there are $n-1$ linearly independent hyperplanes  $(x_i=0)$ for
  $i=1,\dots, n-1$, then $\phi$ is an automorphism of
  $\gm^{n-1}\times \a^1\cong \a^n\setminus(x_1\cdots x_{n-1}=0)$.
  Since any polynomial map $\a^1\to \gm$ is constant, we see that
  $\phi$ is of the form
  $$
  (x_1,\dots, x_n)\mapsto
  \bigl(c_1x_1,\dots, c_{n-1}x_{n-1}, c_nx_n+p(x_1,\dots, x_{n-1})\bigr),
  $$
where $c_i\in k^\times$ and $p\in k[x_1,\dots, x_{n-1}]$.
  \end{exmp}

By contrast, many automorphisms are the identity on a  codimension $\geq 2$  subset.

\begin{exmp}  Let $Z\subset \a^n$ be a closed subset of dimension $\leq n-2$.
  Let 
  $$
  \pi: (x_1, \dots, x_n)\mapsto (x_1, \dots, x_{n-1})
  $$
  be the projection and $W\subset \a^{n-1}(x_1, \dots, x_{n-1})$ the closure of the image of $Z$. Let $p(x_1, \dots, x_{n-1})$ be a polynomial that vanishes on $W$. Then
$$
  \pi: (x_1, \dots, x_n)\mapsto \bigl(x_1, \dots, x_{n-1},x_n+p(x_1, \dots, x_{n-1})\bigr)
  $$  
  is an automorphism of $\a^n$ that is the identity on $Z$. They generate an infinite dimensional subgroup of $\a^n$ that is transitive on $\a^n\setminus Z$.
  \end{exmp}




\section{Lifting of automorphisms}

For the proof of Theorem~\ref{main.k.thm.i}, it is better to use geometric language. Thus
let $\pi: Y\to X$ be morphism of $k$-schemes. We aim to understand which
automorphisms of $X$ lift to automorphisms of $Y$.

\begin{say}[Basic set-up]\label{gen.setup.say}
  Let $\pi: Y\to X$ be a finite, dominant  morphism of  geometrically irreducible and geometrically   normal $k$-schemes.
  Let $Y^{\circ}\subset Y$, $X^{\circ}\subset X$ be  nonsingular, dense, open subsets such that the restriction
  $\pi^\circ: Y^\circ\to X^\circ$ is  a finite and   \'etale morphism.
  
  We study situations when one of the following assumptions holds. 
 \begin{enumerate}
 \item $Y^{\circ}_{\bar k}$ is simply connected.
   \item Every \'etale cover of $Y^{\circ}_{\bar k}$ of degree $\leq \deg \pi$ 
is trivial.
\item Every Galois \'etale cover of $Y^{\circ}_{\bar k}$ of degree dividing $\deg \pi$ 
is trivial, and $\pi$ is Galois.
 \end{enumerate}
 Our main interest is in $Y=\a^n_k:=\spec_k k[{\mathbf x}] $ or  $Y=\hat\a^n_k:=\spec_k k[[{\mathbf x}]]$, and  (1--3) are dictated by what holds for $Y^{\circ}:=Y\setminus Z$, where  $Z$ is a closed subset of dimension $\leq n-2$. Then  $Y^{\circ}$ satisfies one of (1--3)  in the following cases.
 \begin{enumerate}\setcounter{enumi}{3}
 \item If   $Y=\a^n_k$ and $\chr k=0$, then (1) holds.
   \item If  $Y=\hat\a^n_k$ and $\chr k$ arbitrary,    then (1) holds.
   \item If  $\chr k=p>0$ and $\deg \pi<p$, then (2) holds.
\item $\chr k=p>0$,  $\pi$ is Galois and $p\nmid \deg\pi$, then (3) holds.
 \end{enumerate}
 To see this, note that  if $Y$ is regular and
$\dim Z\leq \dim Y-2$, then 
 $\pi_1(Y\setminus Z)=\pi_1(Y)$ by
 Nagata's purity theorem  (cf.\ \cite[\href{https://stacks.math.columbia.edu/tag/0BMB}{Tag 0BMB}]{stacks-project}).

 Next $\pi_1(\hat\a^n_k)=\{1\}$ over any separably closed field
 (see, for example, \cite[\href{https://stacks.math.columbia.edu/tag/04GE}{Tag 04GE}]{stacks-project}). 
 The Hurwitz formula shows that
 $\a^1_{\bar k}$ has no \'etale covers that are tamely ramified at infinity.
 Restricting a cover of $\a^n_{\bar k}$ to a general line through the origin gives the rest our claims.
 \end{say}

The following obvious lifting result applies
only to separably closed base fields.

\begin{lem}\label{all.lift.cor}
  Using the notation  of (\ref{gen.setup.say}), assume that
   $Y^\circ$ is simply connected. Then every automorphisms of $X$ lifts to an automorphism of $Y$.
\end{lem}

Proof.  For  $\phi\in \aut(X)$, let $\tilde Y$ denote the normalization of the fiber product of
$\pi:Y\to X$ and $\phi\circ \pi: Y\to X$.
The first projection $\sigma: \tilde Y\to Y$ is \'etale over $Y^\circ$, hence
$\tilde Y$ is a disjoint union of copies of $Y$ since  $Y^\circ$ is simply connected. Each of these irreducible components gives a lifting. \qed
\medskip

\begin{say}[Multiplicity filtration]\label{mult.filt.say}
Let $(A, m)$ be a complete local $k$-algebra. The kernel of
$\aut_k(A)\to \GL(m/m^2)$ is unipotent.
Thus any reducive subgroup of $\aut_k(A)$ is also a subgroup of
$\GL(m/m^2)$.

The problem is that in most cases the image of $\aut_k(A)\to \GL(m/m^2)$ is much smaller than
$\GL(m/m^2)$, and is hard to determine.
For $A=k[[{\mathbf x}]]^G$, 
the generators of the maximal ideal   are the minimal generators of the ring of invariants. In most cases these are not known, we do not even have  good bounds on the dimension of $ m/m^2$.  
We need to find a representation of $\aut\bigl(k[[{\mathbf x}]]^G\bigr)$
that relates more directly to $G\subset \GL_n$. 

We  understand   much better  
the {\it  multiplicity filtration} of $k[[{\mathbf x}]]^G$. This is given as
$$
I_r:=k[[x_1,\dots, x_n]]^G\cap (x_1,\dots, x_n)^r.
\eqno{(\ref{mult.filt.say}.1)}
$$
That is, 
$I_r:=\{f: \mult_0 f\geq r\}$, 
where $\mult_0 f$ is the multiplicity at the orgin, when we view $f$ as an element of  $ k[[x_1,\dots, x_n]]$.
\end{say}


  \begin{cor}\label{mult.pres.cor}
    Let $k$ be a field and $G\subset \GL_n$ a finite, geometrically reduced  subgroup  scheme without pseudo-reflections.
  Then $\aut\bigl(k[[{\mathbf x}]]^G\bigr)$ preserves the multiplicity filtration.
  \end{cor}

Proof. It is enough to prove this over $\bar k$. Then every automorphisms of $k[[{\mathbf x}]]^G$ lifts to automorphisms of $k[[{\mathbf x}]]$ by Lemma~\ref{all.lift.cor}.

On  $ k[[x_1,\dots, x_n]]$ the multiplicity filtration is the same as the
filtration by the powers of the maximal ideal, so preserved by all automorphisms.
Thus the same holds for $k[[{\mathbf x}]]^G$.\qed
\medskip

We can now prove parts of Theorem~\ref{main.k.thm.i}.

\begin{cor} \label{ex.seq.right.cor} 
Let $k$ be a field and $G\subset \GL_n$ a  finite, geometrically reduced  subgroup scheme  without pseudo-reflections. Then 
  $$
  \aut^{\rm graded}\bigl(k[{\mathbf x}]^G\bigr) \cong 
  \bigl({\bf N}(G,\GL_n) /G\bigr)(k).
  \eqno{(\ref{ex.seq.right.cor}.1)}
  $$
\end{cor}

Proof. Let $k^s\supset k$ be a separable closure.
Since $k[[{\mathbf x}]]^G$ is the completion of
    $k[{\mathbf x}]^G$,  every graded automorphism  $\phi$ of
    $k[{\mathbf x}]^G$ has a natural extension to an automorphism of
    $k[[{\mathbf x}]]^G$, hence lifts to  an automorphism  $\bar \phi$ of
$k^s[[{\mathbf x}]]$  by Lemma~\ref{all.lift.cor}  and (\ref{gen.setup.say}.4).
Since $\bar\phi$ preserves the  powers of the maximal ideal,
it induces an automorphism of the corresponding graded ring. Thus 
 we get
 $\tilde \phi \in \aut^{\rm graded} k^s[{\mathbf x}]\cong \GL_n(k^s)$.
 This shows that
 $$
  \aut^{\rm graded}\bigl(k^s[{\mathbf x}]^G\bigr) \cong 
  \bigl({\bf N}(G,\GL_n(k^s)\bigr) /G(k^s).
  \eqno{(\ref{ex.seq.right.cor}.2)}
  $$
  All this is $\gal(k^s/k)$-equivariant, and, by definition.,
  $\bigl({\bf N}(G,\GL_n) /G\bigr)(k)$ is the group of $\gal(k^s/k)$-fixed points of $\bigl({\bf N}(G,\GL_n(k^s)\bigr) /G(k^s)$. \qed

\begin{cor} \label{ex.seq.cl} 
Let $k$ be a field and $G\subset \GL_n$ a  finite, geometrically reduced  subgroup scheme  without pseudo-reflections. Then there is a split exact sequence 
  $$
  1\to \autu\bigl(k[[{\mathbf x}]]^G\bigr) \to  \aut\bigl(k[[{\mathbf x}]]^G\bigr)\to
  \bigl({\bf N}(G,\GL_n) /G\bigr)(k)\to 1.
  \eqno{(\ref{ex.seq.cl}.1)}
  $$
\end{cor}

  Here $\autu\bigl(k[[{\mathbf x}]]^G\bigr)$ denotes the subgroup of those automorphisms that stabilize  the multiplicity filtration  $\{I_r\}$ of (\ref{mult.filt.say}.1) and act trivially on the graded pieces $I_r/I_{r+1}$.
  On the right hand side, $\bigl({\bf N}(G,\GL_n) /G\bigr)(k) $
  is as in (\ref{main.k.thm.i}.3).

  \medskip

  Proof. By Corollary~\ref{mult.pres.cor}, $ \aut\bigl(k[[{\mathbf x}]]^G\bigr)$
  preserves the multiplicity filtration. Thus we get an exact sequence
  $$
  1\to \autu\bigl(k[[{\mathbf x}]]^G\bigr) \to  \aut\bigl(k[[{\mathbf x}]]^G\bigr)\to
  \aut^{\rm graded}\bigl(\oplus_{r\geq 0}I_r/I_{r+1}\bigr).
  \eqno{(\ref{ex.seq.cl}.2)}
  $$
  Next observe that $\oplus_{r\geq 0}I_r/I_{r+1}\cong k[{\mathbf x}]^G$
    and 
    $$
    \aut^{\rm graded}\bigl(k[{\mathbf x}]^G\bigr)\cong \bigl({\bf N}(G,\GL_n) /G\bigr)(k)
    $$
    by Corollary~\ref{ex.seq.right.cor}.
    Finally $k[[{\mathbf x}]]^G$ is the completion of
    $k[{\mathbf x}]^G$, thus every graded automorphism of
    $k[{\mathbf x}]^G$ has a natural extension to an automorphism of
    $k[[{\mathbf x}]]^G$.  This gives the splitting of
    $ \aut\bigl(k[[{\mathbf x}]]^G\bigr)\to
  \bigl({\bf N}(G,\GL_n) /G\bigr)(k)$. \qed

\begin{rem}  Let $(X,x)$ be a klt singularity over a field of characteristic 0. 
  By  \cite{MR4325260}
  there is a unique valuation centered at $x$ that
  minimizes  the normalized volume. The corresponding filtration on $\o_X$ is thus
  invariant under $\aut (X, x)$. This is interesting when
  this filtration is not equivalent to the `powers of the maximal ideal' filtration.
  This happens for most quotient singularities.

 In our case $G$ acts on $\a^n$ and on
the blow-up  $B_0\a^n$.  Let the quotients be
$\pi:Y_G\to X_G$ with exceptional divisor $E\cong \p^{n-1}/G$.
This $E$ gives the normalized  volume minimizer by \cite[7.1.1]{MR4118616}, giving  the multiplicity filtration.

As a  side remark, note that
going from $X_G$ to $Y_G$ might be  a  step in a resolution algorithm for quotient singularities.
However, as noted in \cite[2.29]{k-res}, repeating it  need not give a resolution.
  For example, starting with the $\z_5$-quotient with resolution dual graph
  $$
  2 \ - \ 3
  $$
  the $E\subset Y_G$ is the $(-1)$-curve in
  $$
  3 \ - \ 2  \ - \ 1 \ - \ 5
  $$
  so we see the same singularity on $Y_G$ that we started with.
  \end{rem}

\section{Lifting unipotent  automorphisms}

The previous lifting results work well for separably close fields.
Over general fields not all automorphisms lift, but we  aim  to show that
all unipotent  automorphisms lift.

\begin{prop}\label{lift.if.poss.prop}
  Using the notation and assumptions of (\ref{gen.setup.say}.1--3), assume  that $Y^\circ(k)$ is nonempty, and 
 set $L:=\pi\bigl(Y^\circ(k)\bigr)$.
  For $\phi\in \aut(X^{\circ})$ the following are equivalent.
  \begin{enumerate}
  \item $\phi$ lifts to $\phi_Y\in \aut(Y^{\circ})$.
  \item $\phi(L)=L$.
  \item $\phi(L)\cap L\neq \emptyset$.
  \end{enumerate}
\end{prop}

Proof. (1) $\Rightarrow$ (2) $\Rightarrow$ (3) are clear.

Next let $\tilde Y$ denote the normalization of the fiber product of
$\pi:Y\to X$ and $\phi\circ \pi: Y\to X$.
The first projection $\sigma: \tilde Y\to Y$ is \'etale over $Y^\circ$.
Liftings of $\phi$ are in one-to-one correspondence with
 sections of $\sigma: \tilde Y\to Y$.

Note that  (3) holds  iff there are $y_1, y_2\in Y^\circ(k)$ such that
$\pi(y_1)=\phi\circ \pi(y_2)$. Then $(y_1, y_2)\in \tilde Y(k)$.
Let $W\subset \tilde Y$ be the unique irreducible component containing 
$(y_1, y_2)$. It has a smooth $k$-point, namely $(y_1, y_2)$, so it is geometrically irreducible.
 The assumptions
(\ref{gen.setup.say}.1--3)  all guarantee that 
 $\sigma^\circ:W^\circ\to Y^\circ$ is an isomorphism over $\bar k$, hence an isomorphism.  \qed

 \medskip

 For finitely generated fields,   $Y^\circ(k)\to (Y^\circ/G)(k)$ is rarely  surjective.
 The following version of Proposition~\ref{lift.if.poss.prop} gives many more cases.

\begin{prop}\label{lift.if.poss.prop.ram}
  We use the notation and assumptions of (\ref{gen.setup.say}.1--3), and fix
  $\phi\in \aut(X)\cap \aut(X^{\circ})$.
Let $D\subset X\setminus X^\circ$ be a geometrically  irreducible divisor,  and $E\subset \pi^{-1}(D)$ a  geometrically  irreducible component.
  Assume that $\pi$ is tamely ramified along $E$, and $\phi$ 
  is the identity on $\o_X/\o_X(-2D)$.

  Then there is a unique lifting of $\phi$  to $\phi_Y\in \aut(Y)$ such that
  $\phi_Y$ is the identity on  $\o_Y/\o_Y(-2E)$.
\end{prop}

Proof. 
As in   the proof of Proposition~\ref{lift.if.poss.prop}, let $\tilde Y$ denote the normalization of the fiber product of
$\pi:Y\to X$ and $\phi\circ \pi: Y\to X$.
Set $L:=k(D)$ and $ K:=k(E)$.
We claim that 
$e:\spec K\into Y$ lifts to $\tilde e:\spec K\into \tilde Y$. 

To see this, let $R_X$ (resp.\ $R_Y$) be the completion of the local ring at the generic point of $D$ (resp.\ $E$).  Then  $R_X\cong L[[s]]$,   $R_Y\cong K[[t]]$
and  $\pi^*s=t^r(\mbox{unit})$ for some $r$ not divisible by $\chr L$.
Write the unit as $c(1+tw(t))$ for some power series $w(t)$. We can replace $t$ by $t/\sqrt[r]{1+tw(t)}$ to get the simpler form
$\pi^*s=ct^r$ for some $c\in K$. So $R_Y\cong K[[s]][t]/(ct^r-s)$. 

$\phi^*s$ is a  power series in $s$, which has the form
$\phi^*s=s(1+sv(s))$  for some power series $v(s)\in L[[s]]$.
The rest is a direct algebra computation, see Lemma~\ref{ram.fib.pr.2}.

Let $W\subset T$ be the unique irreducible component containing
$\tilde e (\spec K)$.
It has a nonsingular,  geometrically irreducible point, so it is geometrically irreducible.
 The assumptions
(\ref{gen.setup.say}.1--3)  all guarantee that 
 $\sigma:W\to Y$ is an isomorphism over $\bar k$, hence an isomorphism.
 Our choice in (\ref{ram.fib.pr.2}) gives that 
  $\phi_Y$ is the identity on  $\o_Y/\o_Y(-(r+1)E)$.
 \qed

\begin{lem} \label{ram.fib.pr.2}
  Let $K/L$ be a finite field extension and   $v(s)\in L[[s]]$  a power series. Assume that $\chr L\nmid r$. 
  Let $B$ be the normalization of 
  $$
  K[[s]][x]/(cx^r-s)\otimes_{L[[s]]}K[[s]][y]/(cy^r-s-s^2v(s)).
  $$
  Then there is a quotient map $q:B\to K[[s]][u]/(cu^r-s)$ such that
  $$
  q(x\otimes 1-1\otimes y)\in (su).
  $$
\end{lem}


Proof.  Since $\chr L\nmid r$, $\sqrt[r]{1+sv(s)}$ is a power series.
Setting $y=z\sqrt[r]{1+sv(s)}$, we get that
$$
K[[s]][x]/(cx^r-s)\otimes_{L[[s]]}K[[s]][z]/(cz^r-s)\cong (K\otimes_kK)[[x,z]]/(cx^r-cz^r).
\eqno{(\ref{ram.fib.pr.2}.1)}
  $$
  This has a quotient $K[[x,z]]/(x-z)\cong K[[x]]\cong  K[[s]][u]/(cu^r-s)$, which is still a quotient after normalization. \qed
  \medskip

  {\it Note \ref{ram.fib.pr.2}.2.}
  If $\phi$ 
    is the identity on $D$, then  $\phi^*s=bs(1+sv(s))$  for some  $b\in L$, and then in    (\ref{ram.fib.pr.2}.1) we need to normalize
    $(K\otimes_kK)[[x,z]]/(x^r-bz^r)$. This involves $\sqrt[r]{b}$, and $\spec K\into X$ lifts iff  $\sqrt[r]{b}\in K$.
    So we really need $\phi$ 
  to be  the identity on $\o_X/\o_X(-2D)$, not just on $\o_X/\o_X(-D)$.
    \medskip

  We can now  show that that every automorphism in
  $\autu\bigl(k[{\mathbf x}]^G\bigr)$ lifts uniquely to
  an automorphism in
  $\autu(k[{\mathbf x}])$.

\begin{thm} \label{main.u.k.thm}
 Let $k$ be a field and $G\subset \GL_n$ a  finite, geometrically reduced  subgroup scheme without pseudo-reflections. 
 Then
 \begin{enumerate}
   \item $\autu(\hat\a^n_k/G)=C\bigl(G, \autu(\hat\a^n_k)\bigr)$, and 
\item
  $\autu(\a^n_k/G)=C\bigl(G, \autu(\a^n_k)\bigr)$ if  $\chr k\nmid |G|$.
  \end{enumerate}
\end{thm}

\medskip

Proof. We discuss (2). The same argument applies to  $\hat\a^n_k$, here the assumption $\chr k\nmid |G|$ is not necessary since  (\ref{gen.setup.say}.5)  is stronger than (\ref{gen.setup.say}.4).

Let $\pi:B_0\a^n_k\to \a^n_k$ denote the blow-up of the origin with exceptional divisor $ E\cong \p^{n-1}$. The $G$-action lifts to $B_0\a^n_k$.
Taking the quotient we get 
$\pi_G:B_0\a^n_k/G\to \a^n_k/G$.

$B_0\a^n_k/G$ is also obtained by blowing up the multiplicity filtration $I_r$ of (\ref{mult.filt.say}.1), which is
invariant under $\aut(\a^n_k/G)$ by Corollary~\ref{mult.pres.cor}.
Therefore every automorphism    $\phi\in \autu(\a^n_k/G)$
lifts to an  automorphism   $\phi_B\in \aut(B_0\a^n_k/G)$ that acts trivially on $E/G$.

We need to show that $\phi_B$ lifts to $\tilde \phi_B\in \autu(\a^n_k)$.
The easy case is when $G$ does not contain scalars.
In this case the $G$-action on $E$ is faithful, so
$\pi_G$ is also \'etale in codimension 1.

If, in addition,  $k$ is infinite, then there is a point  $p\in E\cong \p^{n-1}$  where $\pi_G$ is  \'etale.
then $\phi_B\bigl(\pi_G(p)\bigr)=\pi_G(p)$, so the lifting is
guaranteed by Proposition~\ref{lift.if.poss.prop}.

In general we aim to apply Proposition~\ref{lift.if.poss.prop.ram} with $D=E/G $.  For this we need to know that  $\phi_B$ is trivial on the first order neighborhood of $D$. This is discussed next.
\qed

\begin{say}\label{unip.lift.unip.say}
More generally, 
let $X=\spec R$ be an integral  scheme and  $x\in X$ a closed point with maximal ideal sheaf $m$.
Let  ${\mathcal I}:=\{R=:I_0\supset I_{1}\supset \cdots\}$ be a sub-multiplicative sequence of $m$-primary ideals, that is, $I_a\cdot I_c\subset I_{a+c}$. The blow-up is $B_{\mathcal I}X:=\proj_X (\oplus_j t^jI_j)$, where $t$ is a formal variable and $\deg (t^jI_j)=j$.

The corresponding  graded ring  is $\gr_{\mathcal I}R:=\oplus_j t^j(I_j/I_{j+1})$. Its proj, denoted by $E_{\mathcal I}$, is the exceptional divisor of $B_{\mathcal I}X\to X$.

Assume now that $\phi\in \aut(X)$ stabilizes ${\mathcal I}$,  acts trivially on $\gr_{\mathcal I}R$, and the latter is an integral domain.

  Let $r$ be the smallest integer such that there is an $x\in I_r\setminus I_{r+1}$.  Then   $t^rx\in t^rI_r$ is not a zerodivisor on  $E_{\mathcal I}$.
The   $(t^rx\neq 0)$ chart of the blow-up can be identified with the spectrum of
$S:=R\bigl[ \cup_j x^{-j}I_{rj}\bigr]$, viewed as a subring of the fraction field of $R$. In particular,  $mS=xS$ is a principal ideal.
That is, $(x=0)$ is a defining equation of the exceptional divisor in this chart.

Since the $\phi$-action is trivial on $xR/I_{r+1}$, it is trivial on $xR/(xm)$, hence the lifted
$\phi$-action is trivial on $xS/(xmS)$, which is the restriction of
$\o_{B_{\mathcal I}X}(-D)/\o_{B_{\mathcal I}X}(-2D)$ to our chart.
\end{say}

 The next example shows that the assumption $p\nmid |G|$ in (\ref{main.u.k.thm}.2)---and hence in (\ref{gen.setup.say}.3)---is necessary.

\begin{exmp}  Let   $\pi:Y\to X$ be the
   Artin-Schreier cover of $X:=\a^n(x_1,\dots, x_n)$ given by  $y^p-y-x_1=0$. Then $Y=\a^n(y, x_2,\dots, x_n)$ and 
   $\pi$ is the quotient by the free $\z/p$ action
   $(y, x_2,\dots, x_n)\mapsto (y+1, x_2,\dots, x_n)$.

   Two Artin-Schreier field extensions  $K(y)/(y^p-y-a_i)$ are isomorphic iff
   $a_2-a_1=b^p-b$ for some $b\in K$. This shows that the covers given by 
   $y^p-y-x_1=0$ and $y^p-y-cx_1=0$ are not isomorphic  if $c\neq 1$.
   Thus the automorphisms $(x_1,x_2,\dots, x_n)\mapsto (cx_1,x_2,\dots, x_n)$
   do not lift to $Y$.
   \end{exmp}


\begin{say}[Proof of Theorem~\ref{main.k.thm.i}]\label{main.k.thm.i.pf}
For $k[[{\mathbf x}]]$, Corollary~\ref{ex.seq.cl} gives the split exact sequence 
$$
  1\to \autu\bigl(k[[{\mathbf x}]]^G\bigr) \to  \aut\bigl(k[[{\mathbf x}]]^G\bigr)\to
  \bigl({\bf N}(G,\GL_n) /G\bigr)(k)\to 1.
  \eqno{(\ref{main.k.thm.i.pf}.1)}
  $$
  Next Theorem~\ref{main.u.k.thm}.1 gives 
  $ \autu\bigl(k[[{\mathbf x}]]^G\bigr)\cong  C\bigl(G, \autu(k[[{\mathbf x}]])\bigr) $.
  These together give (\ref{main.k.thm.i}.1).

  For $k[{\mathbf x}]$,  first we observe that
  every automorhism of
  $k[x_1,\dots, x_n]^G$ stabilizes the ideal $(x_1,\dots, x_n)\cap k[x_1,\dots, x_n]^G $ by Lemma~\ref{origin.fixed.lem}, hence  extends to an
  automorphism of $k[[{\mathbf x}]]^G$. Thus (\ref{main.k.thm.i.pf}.1)
  gives
   the split exact sequence 
$$
  1\to \autu\bigl(k[{\mathbf x}]^G\bigr) \to  \aut\bigl(k[{\mathbf x}]^G\bigr)\to
  \bigl({\bf N}(G,\GL_n) /G\bigr)(k)\to 1.
  \eqno{(\ref{main.k.thm.i.pf}.2)}
  $$
  The required identification of $\autu\bigl(k[{\mathbf x}]^G\bigr)$
  is now given by Theorem~\ref{main.u.k.thm}.2, which needs the assumption $\chr k\nmid |G|$. \qed
  \end{say}

\begin{lem} \label{origin.fixed.lem}
  Let $k$ be a field and $G\subset \GL_n$ a  finite, geometrically reduced subgroup scheme  without pseudo-reflections. Assume that  $0\in k^n$ is the only $G$-invariant vector. Then every automorhism of
  $k[x_1,\dots, x_n]^G$ stabilizes the ideal $(x_1,\dots, x_n)\cap k[x_1,\dots, x_n]^G $.
\end{lem}

Proof. It is enough to prove this when $k$ is algebraically closed.
Let $X:=\spec_k k[x_1,\dots, x_n]^G$, $p\in X$ a closed point and
$\bar p\in \a^n$ a preimage of $p$. Let $\hat X_p$ denote the completion of $X$ at $p$. 
Then the fundamental group of
$\hat X_p\setminus \sing X$ is the stabilizer of $\bar p$.
By assumption, the origin is the only point whose stabilizer is $G$. \qed

\section{Normalizers of finite subgroups of $\aut(\a^n_k)$}

Fix a field $k$ and 
let  $G\subset \GL_n$ be a   geometrically reduced subgroup scheme,
where now we view $\GL_n$ as an algebraic group over $k$ that acts linearly on affine $n$-space $\a^n_k$. 
We aim to understand when the centralizer $C\bigl(G, \aut(\a^n_k)\bigr)$ is  small.


The following result, modeled on \cite{MR1288046, MR2340971},  is somewhat complicated, but turns out to apply in many cases.

\begin{thm}\label{bir.ired.act.thm}
    Let $G\subset \GL_n$ be a geometrically reduced $k$-subgroup scheme and
  $\phi\in \aut(\a^n_k)$ an automorphism that commutes with $G$. Then
  \begin{enumerate}
  \item  either $\phi$ extends to an automorphism of $\p^{n}_k$, hence
    $\phi\in \GL_n(k)$,
  \item or there are $G$-equivariant rational maps
    $$
    (G\curvearrowright \p^{n-1}_k) \stackrel{\tau}{\sim}
    \bigl(G\curvearrowright (\p^r_k\times Z)\bigr)  \stackrel{\pi}{\to}
    (G\curvearrowright Z) \stackrel{\sigma}{\to}
    (G\curvearrowright \p^{n-1}_k),
    \eqno{(\ref{bir.ired.act.thm}.2.{\rm a})}
    $$
    where $\tau$ is birational, $r\geq 1$,  and $\pi$ is the coordinate projection. 
  \end{enumerate}
 \noindent  Note that the $G$-action on $ \p^r_k\times Z$ may be only rational, and
    the first projection $\p^r_k\times Z\to \p^r_k$ need not be $G$-equivariant.
\end{thm}

Proof. $\phi$ extends to a birational map $\bar\phi:\p^n_s\map \p^n_t$;
we use the subscripts source/target to distinguish them. Let
$W$ be the closed graph of   $\bar\phi$, with coordinate projections  $p_s, p_t:W\to \p^n$.
Thus $\bar\phi=p_t\circ p_s^{-1}$. Let $H_s\subset \p^n_s$ be the hyperplane at infinity.

If the restriction of $p_t$ to $(p_s^{-1})_*(H_s)\to H_t$ is birational, then
$\bar\phi$ is a local isomorphism at the generic points of $H_s$ and $H_t$. Then
$\bar\phi$ is an  isomorphism by \cite{mats-mumf}; see also
\cite[11.39]{k-modbook} for a more general version.

Otherwise, by a result of Zariski and Abhyankar (see for example \cite[2.45]{km-book}) there is a sequence of $G$-equivariant blow-ups
$$
X_m\to X_{m-1}\to\cdots \to X_0=\p^n_t,
$$
such that the birational transform of $H_s\cong \p^{n-1}_k$ on $X_m$ is the exceptional divisor of the last  blow-up $X_m\to X_{m-1}$. Let $Z\subset X_{m-1}$ be the center of this blow-up. Then we get (2).
\qed
\medskip

{\it Complement \ref{bir.ired.act.thm}.3.} 
  The proof shows that in (\ref{bir.ired.act.thm}.2.{\rm a}) the kernel of
  $G\to \bir(Z)$ acts linearly on $\p^r_{k(Z)}$. 

\medskip

We turn (\ref{bir.ired.act.thm}.2.{\rm a}) into a definition.

\begin{defn}\label{bir.red.defn}
  Let $k$ be a field and  $G\curvearrowright \p^{m}_k$  a group scheme acting on $\p^{m}_k$.  For $m\geq 2$ the action is {\it birationally reducible}
  if there are $G$-equivariant rational maps
    $$
    (G\curvearrowright \p^{m}_k) \stackrel{\tau}{\sim}
    \bigl(G\curvearrowright (\p^r_k\times Z)\bigr)  \stackrel{\pi}{\to}
    (G\curvearrowright Z),
    \eqno{(\ref{bir.red.defn}.1)}
    $$
    where $\tau$ is birational, $0<r<m$,  and $\pi$ is the coordinate projection. Note that the $G$-action on $ \p^r_k\times Z$ may be only rational.

    Otherwise $G\curvearrowright \p^{m}_k$ is  {\it birationally irreducible.}
    Clearly, if $H\subset G$ is a subgroup, and
    $H\curvearrowright \p^{m}_k$ is   birationally irreducible, then so is
    $G\curvearrowright \p^{m}_k$.

     For $G\subset \GL_{m+1}$, let $G\curvearrowright \p^{m}_k$ be the induced action.  If  the $G$-action has an invariant subspace of codimension $\geq 2$, then
   $G\curvearrowright \p^{m}_k$ is birationally reducible   by (\ref{bir.red.exmps}.1).

     If $m=1$ then $0<r<m$ is impossible. So we call $G\curvearrowright \p^{1}_k$ {\it birationally reducible} iff it has a fixed point.
     (This is done mainly to make Corollary~\ref{bir.ired.act.cor} work out.)

  {\it Comment.} This definition suits our current needs, but there are  other versions that may be more natural.

  Requiring $Z\cong \p^{m-r}$ would be a closer analog of the notion of reducible linear representation. In general
  $\p^{m}_k\sim \p^r_k\times Z$ does not imply that $Z\sim \p^{m-r}_k$ \cite{MR86m:14009}, but
  I do not know such  examples that arise from  automorphisms.
\end{defn}

Putting together the definition with Theorem~\ref{bir.ired.act.thm} gives the following.

\begin{cor}\label{bir.ired.act.cor}
  Let $G\subset \GL_n$ be a  geometrically reduced $k$-subgroup scheme whose induced action on $\p^{n-1}_k$ is birationally irreducible.
  Then \begin{enumerate}
    \item $ N\bigl(G, \autu(\a^n_k)\bigr)=C\bigl(G, \autu(\a^n_k)\bigr)=\{1\}$, and 
   \item $ N\bigl(G, \aut(\a^n_k)\bigr)=N\bigl(G, \GL_n\bigr)$.
     \end{enumerate}
\end{cor}

Proof. Let $\phi\in C\bigl(G, \aut(\a^n_k)\bigr)$. We need to show that
$\phi\in \GL_n$.  If this is not the case, then, by Theorem~\ref{bir.ired.act.thm},  we get
$G$-equivariant rational maps as in 
(\ref{bir.ired.act.thm}.2.{\rm a}).
Here $r<n-1$ is impossible since $G\curvearrowright \p^{n-1}_k$ is 
birationally irreducible. Thus $r=n-1$, $Z$ is a point and
$\sigma(Z)\in \p^{n-1}_k$ is a $G$-fixed point.

For $n\geq 3$,  a birationally irreducible action has no fixed points by (\ref{bir.red.exmps}.1). For $n=2$  there are no fixed points by definition.
 \qed
\medskip

Birational irreducibility is only a sufficient condition, but, by Example~\ref{reducible.exmp.1},  linear irreducibility is  almost necessary for Corollary~\ref{bir.ired.act.cor}.

Next we give  examples of actions that are  birationally reducible
or are known to be  birationally irreducible.

\begin{exmp}\label{bir.red.exmps} The following are birationally reducible actions on $\p^n_k$. Here it is more natural to use the classical convention that
  $\p(V)$ is the projective space of all lines through the origin of $V$.

  (\ref{bir.red.exmps}.1) Assume that $G\subset \GL(V)$ is reducible with invariant subspace $W$ of codimension $\geq 2$. Projecting $\p(V)$ from $\p(W)$ gives that
 $$
    (G\curvearrowright \p(V)) \stackrel{\tau}{\sim}
    \bigl(G\curvearrowright \bigl(\p(W+k)\times \p(V/W)\bigr)\bigr).
    $$
    If  $W$ has codimension $1$ and the characteristic is 0, then it has a $1$-dimensional $G$-invariant complement, so $G\curvearrowright \p(V)$ is
    birationally reducible once $\dim V\geq 3$. 
    
    See  Proposition~\ref{pgl.q.rigid.prop} for a   birationally irreducible action with a  codimension 1 invariant subspace in positive characteristic.

    (\ref{bir.red.exmps}.2) Assume that $G\subset \GL(V)$ is irreducible but imprimitive, thus it preserves a direct sum decomposition $V=\oplus_{i=1}^r V_i$.
    Assume that  $\dim V_i>1$.  We get $\p(V_i)\subset \p(V)$.
    Note that a general point  $p\in \p(V)$ is contained in a unique linear space  $L(p)$ of dimension $r-1$ that meets each $\p(V_i)$. Thus we get that
$$
    (G\curvearrowright \p(V)) \stackrel{\tau}{\sim}
    \bigl(G\curvearrowright \bigl(\p(k^r)\times \tprod_i\p(V_i)\bigr)\bigr).
    $$
    See (\ref{no.hom.exmps}.5) for the $\dim V_i=1$ case.
    
    (\ref{bir.red.exmps}.3) With a little more care, the previous construction also works if the $V_i$ are conjugate over $k$.
    We need to replace the factor $\tprod_i\p(V_i)$ with the Weil restriction of $\p(V_1)$.

     (\ref{bir.red.exmps}.4) The  permutation representation of  $S_{n+2}$ and  $A_{n+2}$ in $\GL_{n+1}$ is irreducible and primitive, but  birationally reducible. 

    More generally, assume that $G\curvearrowright \p^n_k$  has an orbit $W$ of $n+2$ points in linearly general position. Let $p\in \p^n$ be a general point. Then there is a unique rational normal curve (of degree $n$) passing through $W\cup\{p\}$ (see, for example, \cite[1.18]{Harris95}).
    This gives a $G$-equivariant morphism
    $\p^n_k\map Z$, whose fibers are the
    rational normal curves passing through $W$. 
    The $n=3$ case is worked out in  \cite[4.7]{MR2605172}.

    If $n$ is odd, or if one of the $n+2$ points is a $k$-point,
then $\p^n_k\map Z$ is a product birationally,  giving 
$$
     (G\curvearrowright \p^{n}_k) \stackrel{\tau}{\sim}
     \bigl(G\curvearrowright (\p^{1}_k \times Z)\bigr).
          $$
    If $n$ is even, then $\p^n_k\map Z$ need not be a product birationally. For $n=2$ there are some birationally irreducible cases; this can be proved as in Proposition~\ref{3d.irred.thm.1}.
    Most likely there are birationally irreducible cases for all even values of $n$. 
\end{exmp}

\begin{exmp}\label{no.hom.exmps} The following are birationally irreducible actions on $\p^n_k$.

  (\ref{no.hom.exmps}.1)  By definition $G\curvearrowright \p^1_k$ is birationally irreducible iff $G$ has no fixed points  over $k$.

  (\ref{no.hom.exmps}.2) If $G\curvearrowright \p^n_k$ is  birationally rigid, then  it is also birationally irreducible; see Definition~\ref{bir.rig.defn}.

  (\ref{no.hom.exmps}.3)
    By   \cite{MR3993279}, $G\subset \PGL_3(\c)$ is 
    birationally rigid  iff $G$ is irreducible and $G\not\cong A_4$, $G\not\cong S_4$.

    (\ref{no.hom.exmps}.4) By  \cite{MR4036497}, $G\subset \PGL_4(\c)$ is 
    birationally rigid  iff $G$ is irreducible, primitive, and $G\not\cong A_5$, $G\not\cong S_5$.

  (\ref{no.hom.exmps}.5) There are many other cases that are not  $G$-birationally rigid, but still birationally irreducible. See
  \cite{https://doi.org/10.48550/arxiv.2202.09319} for such subgroups of $\PGL_4(\c)$.

 (\ref{no.hom.exmps}.6) 
  Let $G$ be a  group acting on  $\p^n_k$ that can not act on a
  rationally connected $(n-1)$-fold over $k$. Then $G\curvearrowright \p^n_k$ is birationally irreducible.

  In positive characteristic, this applies to $\PSL_n(\f_q)$ and
  $\PSU_n(\f_q)$ for $n=3,4$, with finitely many exceptions.
  
  The orthogonal groups $\PSO_n(\f_q)$ act on a quadric hypersurface in $\p^{n-1}$ by definition, so the above argument does not apply.  See Corollary~\ref{pso.cor} for  birationally irreducible cases.

  For $k=\c$ and $n=4$ we have actions of  $\PSp_4(\f_3)$ and  $\PSL_2(\f_{11})$, see \cite{MR1835851}.
  However, these groups both have invariant  Fano hypersurfaces, see  \cite{MR2914804}. The examples in \cite[3.6]{MR3519762} are possible candidates.
  However, \cite{MR3519762} suggests that there may be only few such examples in high dimensions.

   \end{exmp}

More  birationally irreducible examples are obtained using
Lemmas~\ref{k-sz-fix.lem}--\ref{not.faithful.lem}.

\begin{prop}\label{k-sz-fix.lem.c1}
   $\PSL_n(\f_q)\curvearrowright \p^{n-1}_{\bar\f_q}$ is  
  birationally irreducible if $q\neq 2$.
\end{prop}

Proof.  Assume that we have
$\PSL_n(\f_q) \curvearrowright (\p^r_{\bar\f_q}\times Z)$.
For some prime $\ell\neq \chr \f_q$, the maximal torus gives a subgroup  $H\cong (\mu_\ell)^{n-1}$
to which Lemma~\ref{not.faithful.lem} applies. Since $\PSL_n(\f_q)$ is simple,
its action on $Z$ is either trivial or faithful. The latter is impossible by
Lemma~\ref{not.faithful.lem}. Then, by (\ref{bir.ired.act.thm}.3)  we get a faithful, linear $\PSL_n(\f_q)$-action on
$\p^r_{\bar\f_q}$ for $r<n-1$, which is again impossible. \qed

\begin{prop}\label{pso.cor}
  Let $q$ be the power of an odd prime and $n\geq 5$.
  Then
\begin{enumerate}
\item   $\PO_n(\f_q)\curvearrowright \p^{n-1}$ is birationally irreducible.
\item  If $n$ is odd then
  $\PSO^+_n(\f_q)\curvearrowright \p^{n-1}$ is also birationally irreducible.
\end{enumerate}
Here 
  $\PSO^+_n(\f_q)\subset \PSO_n(\f_q)$ denotes the
subgroup of elements with spinor norm 1. It is a simple subgroup of
   index 2.
\end{prop}

Proof. The diagonal group  $\mu_2^n$ is in $\OO_n(\f_q)$.
 So $\mu_2^{n-1}$ is a subgroup of $\PO_n(\f_q)$, and 
 we get the first claim  as in the proof of Proposition~\ref{k-sz-fix.lem.c1}.

 If $n$ is odd then $\PSO_n(\f_q)=\PO_n(\f_q)$, and 
 the spinor norm agrees with the determinant on  $\mu_2^n$.  So $\mu_2^{n-1}$ is a subgroup of $\PSO^+_n(\f_q)$, and
the rest goes as before. \qed

\begin{rem} Note that $\SL_n$ and $\OO_n$ contain many pseudo-reflections, so
  the triviality of $ C\bigl(\SL_n, \autu(\a^n)\bigr)$  and of
  $ C\bigl(\OO_n, \autu(\a^n)\bigr)$ also follows from (\ref{man.ps.r.say}).

  However,  $\SO^+$ does not contain pseudo-reflections, so  we do need
  Proposition~\ref{pso.cor} and Corollary~\ref{bir.ired.act.cor} together
  to see that $ C\bigl(\SO^+_n, \autu(\a^n)\bigr)$ is trivial.
\end{rem}

\medskip

To get another set of examples, let  $T\subset \GL_n(\c)$  be the diagonal subgroup and $N(T)$ its normalizer. Thus $N(T)=T\rtimes S_n$, where we realize $S_n$ as the group of permutation matrices.  For $G\subset S_n$ and $m\in\n$, set
  $G_m:=\mu_m^n\rtimes G\subset N(T)$, and
  let $\bar G_m$ be the image of $G_m$ in $\PGL_n(\c)$.
  We have a split  extension
  $$
  1\to \mu_m^{n-1}\to \bar G_m\to G\to 1.
  $$

\begin{cor}\label{k-sz-fix.lem.c2}
  $\bar G_m\curvearrowright \p^{n-1}_{\c}$ is  birationally irreducible if
  $m>1$ and $G\curvearrowright\mu_{p}^{n-1}$ is irreducible for some prime $p\mid m$.
\end{cor}

Proof.  As in the proof of Corollary~\ref{k-sz-fix.lem.c1} we obtain that
the $\mu_{p}^{n-1}$-action on $Z$ is not faithful.
By our assumptions,  $\mu_{p}^{n-1}\subset \bar G_{p} $ is a minimal normal subgroup.
   Thus the $\mu_{p}^{n-1}$-action on $Z$ is trivial.
   This gives a faithful $\mu_{p}^{n-1}$-action on $\p^r_K$, which is impossible since $r<n-1$. \qed
   \medskip

   {\it Remark \ref{k-sz-fix.lem.c2}.1.}
$G\curvearrowright\mu_{p}^{n-1}$ is irreducible for $p\gg 1$ iff the
$G$-action on the quotient $\c^n/\c(1,\dots, 1)$  is irreducible.

For $n\geq 4$, the action $A_n\curvearrowright\mu_{p}^{n-1}$ is irreducible iff $p\nmid n$.

$\z/n\curvearrowright\mu_{p}^{n-1}$ is irreducible iff
$(x^n-1)/(x-1)$ is irreducible over $\f_p$. That is, 
$n$ is prime, $p\neq n$ and $n\nmid  p-1$.
\medskip

   {\it Remark \ref{k-sz-fix.lem.c2}.2.} The corollary also holds in positive characteristic, as long as $p\neq \chr k$.

\medskip

   {\it Remark \ref{k-sz-fix.lem.c2}.3.}
   The previous method shows that if  $g:\p^{n-1}_K\map Z$ is any
   $\bar G_m$-equivariant morphism to  a smooth variety of dimension $<n-1$, then
   $\mu_{p}^{n-1}$ acts faithfully on the fibers of $g$.
   We have fixed points in some fibers, but these may be singular fibers.

   However, if $p\geq n$ and the generic fibers are separably rationally connected, then  \cite{MR3905117} guarantees fixed points in smooth  fibers, leading to a contradiction as before.

   If $\chr K=0$, we get that 
   $\bar G_{p}\curvearrowright \p^{n-1}_K$ is  not $\bar G_{p}$-birational to any nontrivial  Mori fiber space.

   Nonetheless, $\bar G_{p}\curvearrowright \p^{n-1}_K$ is  most likely not
   birationally rigid. For $n=4$ this is proved in
   \cite[2.1]{MR4036497}.

   \medskip

   The following definition of birational rigidity avoids the definition of Mori fiber spaces.  It is equivalent to the usual definition, whenever
resolution of singularities  and the minimal model  program are known to work.
See \cite{MR1798984} for a more general introduction.

   \begin{defn}\label{bir.rig.defn}[Birational rigidity]     
  A group action  $G\curvearrowright X$ on a smooth projective variety is
  {\it birationally rigid} if the following holds.

  Let $G\curvearrowright Y$ be any action on a smooth projective variety $Y$, that is $G$-birationally equivalent to
  $G\curvearrowright X$. Then, whenever we run the $G$-equivariant minimal model  program on $Y$, the end result $G\curvearrowright Y^{\rm m}$ is isomorphic to $G\curvearrowright X$. (See \cite{km-book} for an introduction to the minimal model  program.)

  It is easy to see that if $G\curvearrowright \p^n_k$ is birationally reducible, then the end result of a $G$-equivariant minimal model program can be a (nontrivial) product $\p^r\times Z'$. So birational rigidity implies birational irreducibility.

  In \cite{MR4036497, MR4036497}, the $G$-birational non-rigidity is established by exhibiting one $G\curvearrowright Y^{\rm m}$ that is not isomorphic to
  $G\curvearrowright \p^n_k$. In most cases, this $Y^{\rm m}$ is not a product.
  However, in many cases the papers do not describe all possible 
  $G\curvearrowright Y^{\rm m}$. So the  birational (ir)reducibility of several of the
  non-rigid cases remains to be decided.
\end{defn}

We already noted the birationally rigid examples from 
\cite{MR3993279, MR4036497}.
Next we discuss an example in positive characteristic, that is
reducible but birationally rigid.

\begin{prop}   \label{pgl.q.rigid.prop}
  Let $P_n(\f_q)\subset \PGL_n(\f_q)$ be the parabolic subgroup that stabilizes the hyperplane $(x_1=0)$.
  The natural action $P_n(\f_q)\curvearrowright \p^{n-1}$ is birationally rigid (over $\bar \f_q$) if $q>n$.
\end{prop}

Proof.  Using \cite[5.12]{ksc}, 
it is enough to prove that,   for any $P_n(\f_q)$-invariant  pencil of hypersurfaces $|H|$ (without fixed components), and  every point $p\in \p^{n-1}$,  we have
$\mult_p |H|< (\deg |H|)/n$.  (Note that the statement in \cite[5.12]{ksc}  uses not the degree of $H$, but
the degree of $H$ divided by the degree of the anticanonical class, which is $n$ for $ \p^{n-1}$.)

We show in (\ref{pgl.q.rigid.prop}.2) that, for every point $p\in \p^{n-1}\setminus (x_1=0)$ and every $P_n(\f_q)$-invariant hypersurfaces $H$, we have $\mult_p H< (\deg H)/q$.

If $p\in (x_1=0)$, then a general member  $H\in |H|$ does not contain
$(x_1=0)$, and, by induction on the dimension,
$\mult_p |H|\leq \mult_p (H|_{(x_1=0)})< (\deg |H|)/q$, and we are done.

\medskip
{\it Claim \ref{pgl.q.rigid.prop}.1.} For any polynomial $f\in \bar \f_q[x_1,\dots, x_n]$
$$
\tsum_{p\in \f_q^n}\ \mult_p f\leq q^{n-1}\deg f.
  $$

  Proof. The base points of the linear system
  $|M|=\bigl|\tsum_i \lambda_i(x_i^q-x_i)\bigr|$ are exactly the points in
  $\f_q^n$. Thus
  $$
  \tsum_{p\in \f_q^n}\ \mult_p f\leq \bigl((f=0)\cdot |M|^{n-1}\bigr)=q^{n-1}\deg f.\qed
  $$

  {\it Claim \ref{pgl.q.rigid.prop}.2.}   Let $H\subset \p^{n-1}$ be a
  $P_n(\f_q)$-invariant hypersurface. Then
  $\mult_pH\leq (\deg H)/q$ for every $p\in \p^{n-1}\setminus (x_1=0)$.

  \medskip
  Proof. Writing $p=(c_1{:}\cdots{:}c_n)$, we may assume that $c_1=1$.
  Acting with the matrices in the unipotent radical of  $P_n(\f_q)$, we get the points $$
  \bigl\{ (1: c_2+a_2 :\cdots : c_n+a_n): a_i\in \f_q\bigr\}.
    $$
    If $H=(h=0)$, then
    $$
    q^{n-1}\mult_pH=\tsum_{p\in \f_q^n}\ \mult_p h(1, x_2+c_2, \dots, x_n+c_n)\leq q^{n-2}\deg H.  \hfill\qed
    $$

    {\it Example \ref{pgl.q.rigid.prop}.3.}  The bound
    $\mult_pH\leq (\deg H)/q$ is essentially optimal.
    If $H$ is the product of all linear forms over $\f_q$, then $\deg H=(q^{n-1}-1)/(q-1)$ and $\mult_p H=(q^{n-2}-1)/(q-1)$.

    \medskip
 {\it Remarks \ref{pgl.q.rigid.prop}.4.}
 $P_n(\f_q)$ contains many pseudo-reflections, so
it would be more interesting to get similar results  for $\SO_n(\f_q)$ and $\SU_n(\f_q)$.

 Also, most likely $P_n(\f_q)\curvearrowright \p^{n-1}$ is birationally rigid even if $q\leq n$,  with finitely many exceptions

\medskip

\medskip

A straightforward generalization of the proof of Theorem~\ref{bir.ired.act.thm} gives the following.

\begin{thm}\label{bir.ired.act.thm.var}
  Let $X$ be a smooth, projective variety  and
  $G\subset \auts(X)$  a geometrically reduced subgroup scheme, such that the $G$-equivariant Picard number of $X$ is $1$. Let   $H\subset X$ be a $G$-irreducible divisor, and 
  $\phi\in \aut(X\setminus H)$  an automorphism that commutes with $G$. Then
  \begin{enumerate}
  \item  either $\phi$ extends to an automorphism of $X$,
  \item or there are $G$-equivariant rational maps
    $$
    (G\curvearrowright H) \stackrel{\tau}{\sim}
    \bigl(G\curvearrowright (\p^r_k\times Z)\bigr)  \stackrel{\pi}{\to}
    (G\curvearrowright Z) \stackrel{\sigma}{\to}
    (G\curvearrowright H),
    \eqno{(\ref{bir.ired.act.thm.var}.2.{\rm a})}
    $$
    where $\tau$ is birational, $r\geq 1$,   $\pi$ is the coordinate projection, and  the $G$-action on $ \p^r_k\times Z$ may be only rational. \qed
  \end{enumerate}
\end{thm}

\section{Low dimensions}

We have a complete description of  $\aut\bigl(k[x,y]^G\bigr)$.
\cite{MR3089030} treats (\ref{2d.irred.thm.1}.1) for $k=\c$,  using an argument  to recognize the linear embeddings $\a^1\into \a^2$. It is not clear how to generalize the method of \cite{MR3089030} to higher dimensions.

\begin{prop}  \label{2d.irred.thm.1}
  Let $k$ be a field and $G\subset \GL_2$ a finite,  geometrically reduced  subgroup scheme.   Assume that $\chr k\nmid |G|$. Then we have one of the following alternatives.
  \begin{enumerate}
  \item  $G$ is not  diagonalizable over $k$.  Then  
  ${\bf N}\bigl(G, \aut(\a^2_k)\bigr)= {\bf N}(G,\GL_2)$ is an algebraic group. Its  identity component is a $k$-torus of dimension 2 if $G$ is commutative, and the group of scalars $\gm$ if $G$ is not commutative.
  \item  $G$ is diagonalizable over $k$ with eigenbasis $(1,0), (0,1)$. 
    \begin{enumerate}
    \item If $G\cap (\gm\times\{1\})$ and $G\cap (\{1\}\times\gm)$ are both nontrivial, then $N\bigl(G, \aut(\a^2_k)\bigr)=\gm\times\gm$.
      \item Otherwise $N\bigl(G, \aut(\a^2_k)\bigr)$ is infinite dimensional.
    \end{enumerate}
    \end{enumerate}
\end{prop}

Proof.  
If $G$ is not  diagonalizable over $k$, then the induced action on $\p^1_k$ is
birationally irreducible by  (\ref{no.hom.exmps}.1).
Thus $ C\bigl(G, \autu(k[V])\bigr)=\{1\}$ by  Corollary~\ref{bir.ired.act.cor}, and   we have (\ref{2d.irred.thm.1}.1).

We discussed the $n$-dimensional version of (\ref{2d.irred.thm.1}.2.a) in  Example~\ref{man.ps.r.say}.
Otherwise $G\cong \mu_n$ acting as $(x,y)\mapsto (\epsilon x, \epsilon^q y)$ for some $0\leq q<n$. If $(q,n)=1$ then
$G$ contains no pseudo-reflections, and $N\bigl(G, \autu(\a^2_k)\bigr)$
contains the elementary automorphisms
$$
(x,y)\mapsto \bigl(cx, y+x^qg(x^n)\bigr)\qtq{and}
(x,y)\mapsto \bigl(x+y^{q'}g(y^n), cy\bigr),
$$
where $0<q'<n$ is the inverse of $q$ modulo $n$.
The arguments in \cite{https://doi.org/10.48550/arxiv.2210.12781} should give that, together with the linear   automorphisms, they  generate  the  full automorphism group.

If $(q,n)>1$, then
$G$ contains a pseudo-reflection, and  $N\bigl(G, \autu(\a^2_k)\bigr)$
consists of the automorphisms
$(x,y)\mapsto \bigl(cx, y+x^qg(x^n)\bigr)$ by  Example~\ref{man.ps.r.say}.
\qed

\medskip
Next we consider  irreducible subgroups $G\subset \GL_3$.
The main work is done in  \cite{MR3993279}, which in turn 
relies on the classification of finite subgroups of
$\GL_3(\c)$, due to \cite{blichfeldt}.
(For lists in positive characteristic, see
\cite{MR2151423}.)

\begin{prop}  \label{3d.irred.thm.1}
  Let $G\subset \GL_3(\c)$ be an  irreducible subgroup. Then
  $N\bigl(G, \aut(\c^3)\bigr)= N\bigl(G,\GL_3(\c)\bigr)$.
  
\end{prop}

Proof. By  \cite{MR3993279}  $G\subset \GL_3(\c)$ is birationally rigid iff it is irreducible,    and $G$ is not a central extension of $A_4$ or $S_4$. Since birationally rigid implies birationally irreducible, it remains to deal with the unique irreducible $A_4$ and $S_4$ actions.
(These are obtained as the permutation representation on $\c^4$  restricted to the hyperplane $\tsum x_i=0$.)

The orbits are described in \cite[3.14]{MR3993279}, and   \cite[2.1]{MR3993279} summarizes the possible links, due to   \cite{MR1311348}.

The length 3 orbit gives the quadratic Cremona transformation. Blowing up the points in a length 6 
orbit  gives a cubic surface, the other side of the link contracts the 6 conics that pass through 5 of the 6 points;  we are  back to $\p^2$. Finally, blowing up the points in a length 4 orbit,   we get a conic bundle
$B_4\p^2\to C\cong \p^1$ with 3 singular fibers, where the 6 curves in the  3 singular fibers form a  single $A_4$-orbit.
 We can now do a sequence of
Type II links that preserve the conic bundle structure and the 3 singular fibers. The only other possibility is a Type I link that takes us back to $\p^2$.
Thus the only possibility for (\ref{bir.ired.act.thm.var}.2.{\rm a}) is
$$
(A_4\curvearrowright \p^2) \stackrel{\tau}{\sim}
    \bigl(A_4\curvearrowright B_4\p^2\bigr)  \stackrel{\pi}{\to}
    (A_4\curvearrowright C) \stackrel{\sigma}{\to}
    (A_4\curvearrowright \p^2).
    \eqno{(\ref{3d.irred.thm.1}.1)}
    $$
    Note, however, that the $A_4$-action   on $C$ reduces to a $\z/3$-action.
    This is clear if we write the 4 points as  $(\pm 1{:}0{:}1), (0{:}\pm 1{:}1)$. Then $B_4\p^2\to C\cong \p^1$ is given by the pencil
    $\lambda (x^2+y^2-z^2)+\mu xy$. Every member is invariant under
    $(x{:}y{:}z)\mapsto (-x{:}-y{:}z)$ and $(x{:}y{:}z)\mapsto (y{:}x{:}z)$.

    However,   the action of $(\z/2)^2\subset A_4$ on $\p^2$ is trivial only at 3 points,  so $\sigma$ in (\ref{3d.irred.thm.1}.1) can not exist.
    The same argument applies for $S_4$. \qed

    \medskip

    Over a number field, we give a  3-dimensional, irreducible example.

  \begin{exmp} Let $p=1+3m$ be a prime and set $\epsilon:=e^{2\pi i/p}$.
    The Galois group of $\q(\epsilon)/\q$ is  $\z/(3m)$.
    Choose $c\neq 1$ such that $c^3\equiv 1\mod p$.
    Let  $\sigma_i$ be the elementary symmetric polynomials in
    $\epsilon,\epsilon^c,\epsilon^{c^2}$  (so  $\sigma_3=1$).
        They generate a  subfield  $K_p\subset \q(\epsilon)$
    such that $\deg \q(\epsilon)/K_p=3$, and $\epsilon, \epsilon^c,\epsilon^{c^2}$ satisfy the equation
    $$
    g_p(x):=x^3-\sigma_1x^2+\sigma_2x-\sigma_3=0.
    $$
    Consider the matrix
    $$
    M_p:=
    \begin{pmatrix}
  0 & 0 & \sigma_3\\
  1 & 0 & -\sigma_2\\
    0 & 1 & \sigma_1
    \end{pmatrix}.
$$
    Its characteristic polynomial is $g_p(x)$, hence its eigenvalues are $\epsilon,\epsilon^c,\epsilon^{c^2}$. $M_p$ defines an irreducible  linear representation
    $\z/p\into \GL_3(K_p)$.
    We claim that the induced   $M_p:\z/p\curvearrowright \p^2_{K_p}$
 are birationaly irreducible, but not birationaly rigid.

    The 3 eigenvectors of $M_p$ give a point of degree 3 on $\p^2_{K_p}$, the other
    $\z/p$-orbits of points have degree divisible by $p$.

    Arguing as in  \cite{MR3993279}  (or as in Proposition~\ref{3d.irred.thm.1}), the  degree 3 point gives the quadratic Cremona transformation.
    This corresponds to the complex conjugate $\bar M_p$, whose eigenvalues are the inverses of the eigenvalues of $M_p$.
    The two actions on  $\p^2_{K_p}$ are not isomorphic.
    The only other possibility is for $p=7$, where   degree 7 points give Bertini involutions.

    For $p=7$ we have  $K_p=\q(\sqrt{-7})$ and 
    $$
    g_3(x)=x^3-\alpha x^2-(1+\alpha)x-1,\qtq{where} 
    \alpha:=\epsilon+\epsilon^2+\epsilon^4=\tfrac{-1+\sqrt{-7}}{2}.
    $$
\end{exmp}

The following is a 4-dimensional imprimitive example.

  \begin{exmp}  \label{2i2i2.exmp}
    Let $2\rm{I}\times 2\rm{I}\subset \GL_4(\c)$ be the block-diagonal subgroup and $\rho$ the involution interchanging the 2 blocks:
    $\rho(x_1{:}x_2{:}y_1{:}y_2)= (y_1{:}y_2{:}x_1{:}x_2)$.
    Together they generate $G\subset \GL_4(\c)$ of order  $2\cdot 120^2$.
    The commutator subgroup of $G$ is $2\rm{I}\times 2\rm{I}$, whose
    center is $\mu_2\times \mu_2$.

    The quotient by $2\rm{I}\times 2\rm{I}$ is the product
    $$
    \c^4/(2\rm{I}\times 2\rm{I})\cong (u_1^2+u_2^3+u_3^5=0)\times (v_1^2+v_2^3+v_3^5=0),
    $$
    and then we take the quotient by
    $\bar \rho:(u_1,u_2,u_3)\leftrightarrow (v_1,v_2,v_3)$.

    $G$ is irreducible but imprimitive, so it is  birationally reducible by (\ref{bir.red.exmps}.2). Nonetheless, we claim that
    $$
  \aut \bigl(k[x_1,x_2,y_1,y_2]^G\bigr)=
  \bigl({\bf N}(G,\GL_4) /G\bigr)(k)\cong k^\times.
   \eqno{(\ref{2i2i2.exmp}.1)}
  $$
    To see this, note that the induced action  on $\p^3$ is faithful on $\bar G:=G/\mbox{(diagonal $\mu_2$)}$.
    The center of  $\bar G$ is the involution
    $\chi: (x_1{:}x_2{:}y_1{:}y_2)\mapsto (-x_1{:}-x_2{:}y_1{:}y_2)$.

    The key observation is that while  $\bar G/(\chi)=G/(\mu_2\times \mu_2)$ does act faithfully on $\p^1\times \p^1$, $\bar G$ can not  act faithfully on
    a rational surface  \cite{MR2641179}. Thus if we have
    $$
    (\bar G\curvearrowright \p^{3}) \stackrel{\tau}{\sim}
    \bigl(\bar G\curvearrowright (\p^r\times Z)\bigr)  \stackrel{\pi}{\to}
    (\bar G\curvearrowright Z) \stackrel{\sigma}{\to}
    (\bar G\curvearrowright \p^{3})
    $$
    as in (\ref{bir.ired.act.thm}.2.{\rm a}), then $\chi$ acts trivially on
    $Z$, hence also on $\sigma(Z)\subset \p^3$. However, the fixed locus of $\chi$ on $\p^3$ is the union of 2 disjoint lines that are interchanged by $\rho$. Thus
    there is no $\rho$-equivariant map from  $\p^3$ to the fixed locus of $\chi$.

    By Theorem~\ref{bir.ired.act.thm} we get the first equality in (\ref{2i2i2.exmp}.1). Automorphisms of $G$ fix the  commutator subgroup  $2\rm{I}\times 2\rm{I}$, thus they are all inner. So
    ${\bf N}(G,\GL_4) /G={\bf C}(G,\GL_4) /G=\gm/\mbox{(diagonal $\mu_2$)}$.\qed

    \medskip

    It is possible that the automorphism group of
    $\c^4/(2\rm{I}\times 2\rm{I})$ is the obvious  $(k^\times)^2\rtimes \z/2$, but the above methods do not show this.
    
The same argument should apply to several other imprimitive $G\subset \GL_4(\c)$ with 2 blocks, but one needs to exclude accidental actions on rational surfaces using the lists of \cite{MR2641179}. 

  \end{exmp}

\section{Fixed points and rational maps}

The following lemmas were used in the proofs of
Propositions~\ref{k-sz-fix.lem.c1}--\ref{k-sz-fix.lem.c2}.

\begin{lem}\cite[A.2]{k-sz-fix}\label{k-sz-fix.lem} Let
$A$ be an abelian group  and
  $g:X\map Z$  an
  $A$-equivariant rational map of  proper $k$-varieties.
  If $A$ has a smooth (closed) fixed point on $X$, then
  it also has a (closed) fixed point on $Z$. \qed
\end{lem}

\begin{lem}\label{k-sz-fix.lem.2} Let
  $A$ be an abelian group  acting faithfully on a $k$-variety $Z$ with a 
  smooth (closed) fixed point. Assume that $\chr k\nmid |A|$.
  Then $A$ can be generated by $\dim Z$ elements.
\end{lem}
 
Proof. We may assume that $k$ is separably closed. The $A$-action on the tangent space at the fixed point is also faithful and diagonalizable. \qed

\medskip

While  \cite[A.2]{k-sz-fix} guarantees a fixed point even if $Z$ is singular,
Lemma~\ref{k-sz-fix.lem.2} needs the fixed point to be smooth.
This is usually not a problem if $A$-equivariant resolution of singularities holds for $Z$. When resolution is not known, the following variant is  useful.

\begin{lem}\label{not.faithful.lem}
  Let $k$ be a field, $A:=\mu_\ell^n$  for some prime $\ell\neq \chr k$, and $g:X\map Z$ an $A$-equivariant rational map of $k$-varieties.
  Assume that $A\curvearrowright X$ has a smooth, isolated, fixed point, and
  $\dim Z<n$. Then $A\curvearrowright Z$ is not faithful.
\end{lem}

Proof. We may assume that $k$ is separably closed.
Let $x\in X$ be a smooth, isolated fixed point. Choosing an $A$-eigenbasis of $m_x/m_x^2$, and $A$-equivariantly lifting it, we get
local coordinates $(h_1, \dots, h_n)$ at $x$ such that $A$ acts on them as
$$
(\epsilon_1, \dots, \epsilon_n)\times (h_1, \dots, h_n)\mapsto
(\epsilon_1h_1, \dots, \epsilon_nh_n).
$$
We use induction on $n$; the claim is clear if $n=1$.

We may assume that $Z$ is normal and projective.
Set $X_1:=(h_1=0)$ and let $A\supset A_1\cong \mu_\ell^{n-1}$ be the kernel of the $A$-action on $h_1$.  Then $g$ restricts to an 
$A$-equivariant rational map  $g_1:X_1\map Z$.
If $g_1$ is dominant, then $A\curvearrowright Z$ is not faithful
since $A\curvearrowright X_1$ is not faithful.

If $g_1$ is not dominant, then, by \cite[2.22]{kk-singbook}, after some blow-ups we may assume that the closure of $g_1(X_1)$ is a divisor, call it
$Z_1\subset Z$.

If $\dim Z<n$, then $\dim Z_1<n-1$, hence $A_1\curvearrowright Z_1$ is not faithful by induction. Thus the kernel $K_1$ of $A\curvearrowright Z_1$
is $\mu_\ell^r$ for some $r\geq 2$.

The kernel acts on $\o_Z(-Z_1)/\o_Z(-2Z_1)$ by a character
$K_1\to k^\times$. Thus there is a $\mu_\ell\cong K_2\subset K_1$ that acts
trivially on $\o_Z/\o_Z(-2Z_1)$. Then $K_2$ acts trivially on $Z$. \qed

\begin{rem}  The above lemmas use only the (non)existence of fixed points.
  The theory of  the equivariant Burnside group
  developed in \cite{https://doi.org/10.48550/arxiv.2010.08902, https://doi.org/10.48550/arxiv.2108.00518, https://doi.org/10.48550/arxiv.2204.03106}
  looks at much finer properties, involving the fixed point sets and the
  representations on their normal bundles. 
It gives a powerful  method to distinguish  $G$-actions modulo $G$-birational equivalence. It would be interesting to see what it says about birational irreducibility.
\end{rem}

\section{Twisted forms of quotient singularities}

The arguments to show Corollary~\ref{main.thm} are pretty standard. The Galois cohomology part goes as expected 
following  \cite[\S.III.1]{MR1466966}, and so is computing  automorphisms using the universal cover. For example, it was used in  \cite{MR3089030} for the same purposes.

\begin{say}[Proof of Corollary~\ref{main.thm}]\label{main.thm.pf}

For any $k$-scheme $Y$, its twisted forms are classified by
$H^1\bigl(\gal({\bar k}/k), \aut(Y_{{\bar k}})\bigr)$; see \cite[\S.III.1]{MR1466966}. Thus we need to understand $\aut\bigl(\bar k[[{\mathbf x}]]^G\bigr)$, which is given by the sequence (\ref{main.k.thm.i}.1)
$$
  1\to C\bigl(G, \autu(\bar k[[{\mathbf x}]])\bigr) \to  \aut \bigl(\bar k[[{\mathbf x}]]^G\bigr)\to
  \bigl({\bf N}(G,\GL_n) /G\bigr)(\bar k)\to 1.
  $$
On the right hand side we have
$$
\bigl({\bf N}(G,\GL_n) /G\bigr)(\bar k)=
N\bigl(G, \GL_n(\bar k)\bigr)/G(\bar k),
$$
which is exactly the group of graded auomorphisms of $k[{\mathbf x}]^G$.
This gives the natural map (\ref{main.thm}.1) $\to$ (\ref{main.thm}.2).

For the converse we need to show that
$$
H^1\Bigl(\gal({\bar k}/k),
 C\bigl(G, \autu(k[[{\mathbf x}]])\bigr)\Bigr)=\{1\}.
$$
 Over a perfect field $k$,  a unipotent group has a decreasing filtration whose subquotients are isomorphic to $\ga$, and  $H^1\bigl(\gal(\bar k/k), \ga\bigr)=0$
by \cite[\S.II.1.1]{MR1466966}. Thus 
$H^1\bigl(\gal(\bar k/k), U\bigr)$ of  any unipotent group $U$ is trivial by \cite[\S.I.5.4]{MR1466966}. \qed
\end{say}
\medskip

\begin{exmp} \label{exmp.2} We give 2 examples when there are no twisted forms.

 (\ref{exmp.2}.1)  Assume that $G\subset \GL_n(\bar k)$ acts irreducibly,
  and every automorphism of $G$ is inner. Then $N\bigl(G, \GL_n(\bar k)\bigr)=G\cdot \gm(\bar k)$, where $\gm(\bar k)$ is the group of scalars. Thus $N\bigl(G, \GL_n(\bar k)\bigr)/G$ is the quotient of $\gm$ by the center of $G$, hence isomorphic to $\gm$. Thus 
  $$H^1\bigl(\gal(\bar k/k), N\bigl(G, \GL_n(\bar k)\bigr)\bigr)=H^1\bigl(\gal(\bar k/k), \gm\bigr)=1,$$ and $k[[{\mathbf x}]]^G$ has no twisted forms.

   (\ref{exmp.2}.2)  Assume that $G$ is simple,
   every automorphism of $G$ is inner, and $\chr k=0$.
  Any representation of $G$ decomposes as  $V=\oplus_i \bigl(V_i\otimes_{\bar k} \bar k^{m_i}\bigr)$
  where the $V_i$ are distinct, irreducible representations.
  Then  $N\bigl(G, \GL_n(\bar k)\bigr)\cong \times_i \GL_{m_i}(\bar k)$. Thus
$H^1\bigl(\gal(\bar k/k), N\bigl(G, \GL_n(\bar k)\bigr)\bigr)=\times_i H^1\bigl(\gal(\bar k/k), \GL_{m_i}(\bar k)\bigr)=1$ by \cite[\S.III.1.1]{MR1466966}.
Therefore 
  $k[[{\mathbf x}]]^G$ has no twisted forms for all
  representations  of $G$. 
  \end{exmp}

The following is closely related  to \cite[6.18]{bre-vis-2}.

\begin{exmp} \label{exmp.1}  Let $G\subset \GL_n$ be a subgroup and assume that every automorphism of $G$ is inner. Then $C\bigl(G, \GL_n(\bar k)\bigr)$, the centralizer of $G$ in $\GL_n$, surjects onto  $N\bigl(G, \GL_n(\bar k)\bigr)$. If, in addition,  the center of $G$ is trivial, then
  $C\bigl(G, \GL_n(\bar k)\bigr)\to N\bigl(G, \GL_n(\bar k)\bigr)/G(\bar k)$ is an isomorphism. Thus we have a natural map
  $$
  \begin{array}{l}
  H^1\Bigl(\gal(\bar k/k), N\bigl(G, \GL_n(\bar k)\bigr)/G(\bar k)\Bigr)\cong 
  H^1\Bigl(\gal(\bar k/k), C\bigl(G, \GL_n(\bar k)\bigr)\Bigr)\to \\
  \qquad\qquad \to
  H^1\bigl(\gal(\bar k/k), \GL_n(\bar k)\bigr)=1.
  \end{array}
  $$
  We conclude that every twist of $k[[{\mathbf x}]]^G$ is of the form $k[[{\mathbf x}]]^{G'}$,
  where $G'\subset \GL_n$ is a twisted form of $G$. 
(More generally, the above argument applies whenever $N\bigl(G, \GL_n(\bar k)\bigr)\to N\bigl(G, \GL_n(\bar k)\bigr)/G$ has a canonical splitting.) 

  Applying the  Nishimura lemma  to $\hat \a^n\to \spec_k  R^{G'}$ we get that  every  resolution of every twist of $k[[{\mathbf x}]]^G$ has a $k$-point over the origin.  Thus $\spec k[[{\mathbf x}]]^G$ is an {\it R-singularity} in the terminology of
   \cite[6.11]{bre-vis-2}.

  See \cite[p.183]{ksc}  for an elegant proof of the Nishimura lemma due to E.~Szab\'o.
  \end{exmp}

\begin{rem} For surface singularities, $N\bigl(G, \GL_2(\bar k)\bigr)$ acts on the dual graph of the minimal resolution, and, except for the $\a^2/\tfrac1{n}(1,1)$ quotients,  the kernel is a  torus. Since  $H^1\bigl(\gal(\bar k/k), \gm\bigr)=0$
  by \cite[\S.II.1.1]{MR1466966}, the twisted forms are determined by the symmetries of the dual graph, except when the dual graph consists of a single curve.

  There is no minimal resolution in higher dimensions, and it is not clear whether there is a similarly clear geometric description.
\end{rem}

\section{Uniqueness of cones}

 We end with a  discussion of general cones.

\begin{defn}\label{gen.cone,defn}
  Let $k$ be a field, $X$  a normal, projective  $k$-variety and $\Delta$  an ample $\q$-divisor on $X$. The {\it cone} over   $(X,\Delta)$ is 
  $$
  R(X,\Delta):=\oplus_{m\geq 0} H^0\bigl(X, \o_X(\rdown{m\Delta})\bigr),
  \qtq{or} 
  C_a(X, \Delta):=\spec_k R(X,\Delta).
  \eqno{(\ref{gen.cone,defn}.1)}
  $$
 The  completions are denoted by $\hat R(X,\Delta)$ and  $\hat C_a(X, \Delta)$.

Both $R(X,\Delta)$ and $\hat R(X,\Delta)$ come with a positive grading, where
$\deg r=m$ iff $r\in H^0\bigl(X, \o_X(\rdown{m\Delta})\bigr)$.
A positive grading is equivalent to a postive $\gm$-action
defined on the homogeneous elements by $(\lambda, r)\mapsto \lambda^{\deg r}r$.
 
Set $L:=\o_X(\rdown{\Delta})$ and let  
$\{\Delta\}:=\Delta-\rdown{\Delta}$  be the  fractional part of $\Delta$.
Then
$$
  R(X,\Delta)\cong\oplus_{m\geq 0} H^0\bigl(X, L^{[m]}(\rdown{m\{\Delta}\})\bigr),
  \eqno{(\ref{gen.cone,defn}.2)}
  $$
  thus the cone is determined by $\bigl(X, \o_X(\rdown{\Delta}), \{\Delta\}\bigr)$. 
\end{defn}

  Conversely, an observation going back at least to \cite{dolgachev, pinkham} is the following.

  \begin{lem}\label{cone.str.lem}  Let $(R, m, \deg)$ be a   normal, complete, local $k$-algebra  with a positive grading such that $R/m$ is finite over $k$.
    Then  $(R, m, \deg)$ is a completed cone over a unique $\bigl(X, \o_X(\rdown{\Delta}), \{\Delta\}\bigr)$.
  \end{lem}

  Sketch of proof. 
 $R$ is the completion of $\oplus_{m\geq 0} R_m$ where $R_m\subset R$ is the set of degree $m$ homogeneous elements.
Set $X:=\spec \oplus_{m\geq 0} R_m$. Then
$R_m=H^0\bigl(X, \o_X(m)\bigr)$ and $R$ is the completion of
$\oplus_{m\geq m} H^0\bigl(X, \o_X(m)\bigr)$.

If we work with  $\hat R(X,\Delta)$, then this way we recover
$\rdown{\Delta}\sim \o_X(1)$. Finally we get 
$m\{\Delta\}$  as the divisor associated to 
$\coker\bigl[\o_X(1)^{\otimes m}\to \o_X(m)\bigr]$
for $m$ sufficiently divisible. \qed




\medskip

Next we consider when  $\hat R(X, \Delta)$
admits a unique positive grading.
Background results on
algebraic group actions on complete local rings, and the definition of
the maximal reductive quotient $\auts^{\rm red}_k\bigl(\hat R(X, \Delta)\bigr)$ are discussed in 
(\ref{aut.A.say}).  (Note that $\auts^{\rm red}_k\bigl(\hat R(X, \Delta)\bigr)$ denotes an algebraic group whose $k$-points are $\aut^{\rm red}_k\bigl(\hat R(X, \Delta)\bigr)$. The difference becomes important for us only when $k$ is finite.)

\begin{prop}\label{cones.unique.pr}
   Let $k$ be a perfect  field, $X$  a normal, projective  $k$-variety and $\Delta$  an ample $\q$-divisor on $X$.
  The following are equivalent.
  \begin{enumerate}
    \item  $\hat R(X, \Delta)$ uniquely determines the triple
      $\bigl(X, \o_X(\rdown{\Delta}), \{\Delta\}\bigr)$.
         \item There is no faithful $\gm^2$-action on $\hat R(X, \Delta)$.
       \item $\auts^{\rm red}_k\bigl(\hat R(X, \Delta)\bigr)$ is  a finite group  extended by $\gm$.
  \item There is no faithful $\gm$-action on $\bigl(X, \o_X(\rdown{\Delta}), \{\Delta\}\bigr)$.
\end{enumerate}
\end{prop}

Proof.
Let $\chi_i:\gm\to \auts_k\bigl(\hat R(X, \Delta)\bigr)$ be 2 commuting $\gm$-actions such that
$\chi_1$ is positive. Then $\chi_2\chi_1^r$ is positive for all $r\gg 1$,
giving infinitely many different positive gradings on  $\hat R(X, \Delta)$.
 This shows that (1) $\Rightarrow$ (2).

If (2) holds then 
the  maximal split torus $T$ of $\auts^{\rm red}_k\bigl(\hat R(X, \Delta)\bigr)$ has dimension 1.
Then the identity component of $\auts^{\rm red}_k\bigl(\hat R(X, \Delta)\bigr)$ is either  $T$ or $\SL_2$. However, it can not be 
$\SL_2$ since in every representation of  $\SL_2$ its torus acts with both positive and negative weights. This proves (2)  $\Rightarrow$ (3), and
(3)  $\Rightarrow$ (1) follows from Lemma~\ref{cone.str.lem},


If there is a faithful $\gm$-action on $\bigl(X, \o_X(\rdown{\Delta}), \{\Delta\}\bigr)$, then we get a faithful  $\gm^2$-action on $\hat R(X, \Delta)$, and
conversely. Thus  (2) $\Leftrightarrow$ (4). \qed


\begin{say}\label{aut.A.say}
  Let $(A, m)$ be a complete local $k$-algebra such that $A/m$ is finite over $k$. Then each
  $\auts_k(A/m^r)$ is an algebraic group and we have natural homomorphisms
  $$
  q^{r+s}_r: \auts_k(A/m^{r+s})\to \auts_k(A/m^r)
  $$
  with unipotent  kernels for $r\geq 2$. For  fixed $r$, the images of
  $q^{r+s}_r$ stabilize for $s\gg 1$, let us denote this subgroup by
  $$
  \auts^{\ell}_k(A/m^r)\subset \auts_k(A/m^r).
  $$
  We identify the automorphism scheme $\auts(A, m)$ with the inverse limit diagram
  $$
  \cdots \to \auts^{\ell}_k(A/m^{r+1})\to \auts^{\ell}_k(A/m^r)\to \cdots
  \eqno{(\ref{aut.A.say}.1)}
  $$
       {\it Warning.} Note that $\aut(A, m)$ is the inverse limit of the $k$-points
       $$
  \cdots \to \aut^{\ell}_k(A/m^{r+1})\to \aut^{\ell}_k(A/m^r)\to \cdots
  \eqno{(\ref{aut.A.say}.2)}
  $$
  but we do not discuss the possible scheme structure of     $\auts_k(A, m)$ here.
  \medskip

  By a finite dimensional subgroup scheme  of  $\auts_k(A, m)$ we mean a
  sequence of  subgroup schemes  $G_r\subset \auts^{\ell}_k(A/m^r)$
  such that the  induced $q^{r+s}_r:G_{r+s}\to G_r$ are surjective and
  isomorphisms for $r\gg 1$. Similar definition applies to quotients.

  Let $U^{\ell}_r\subset \aut^{\ell}_k(A/m^r)$ denote the unipotent radical and
  $R^{\ell}_r:=\aut^{\ell}_k(A/m^r)/U^{\ell}_r$ the maximal reductive quotient.
  Since the kernel of $q^{r+s}_r$ is unipotent,
  the induced maps  $R^{\ell}_{r+1}\to R^{\ell}_r$ are isomorphisms for $r\geq 2$.
  This defines  $\auts^{\rm red}_k(A,m):=R^{\ell}_r$ (for any $r\geq 2$) as the
  maximal reductive quotient of  $\auts_k(A, m)$.

  If the characteristic is 0, then  the Levi subgroup realizes $\auts^{\rm red}_k(A,m)$ as a subgroup of 
$\auts_k(A,m)$, unique up to conjugation.  Thus the splitting of (\ref{ex.seq.cl}.1)
  is automatic.
  (In positive  characteristic there are no Levi subgroups in general.)
  
  If $k$ is perfect  then a  Borel subgroup of  $\auts^{\ell}_k(A/m^{r+1})$ maps onto a
   Borel subgroup of  $\auts^{\ell}_k(A/m^{r})$,  and the dimension of their maximal split subtori is the same for $r\geq 2$  \cite[20.9]{borel}. 
   Taking the limit we get maximal split subtori  $T\subset \auts_k(A,m)$.
   They are conjugate to each other and map isomorphically to the
   maximal split subtori  of  $\auts^{\rm red}_k(A,m)$.
  \end{say}
  




   \begin{rem} If there is a $\gm^r$-action on $\bigl(X, \o_X(\rdown{\Delta}), \{\Delta\}\bigr)$, then the above argument shows that we get a
     $\gm^{r+1}$-action on $\hat R(X, \Delta)$, and all other cone structures on $\hat R(X, \Delta)$ are obtained by quotienting by some
     $\gm\into \gm^{r+1}$.

     If $X$ is smooth and not uniruled,
     then $X$ gives the unique non-uniruled exceptional divisor over the vertex of $\hat R(X, \Delta)$,  cf.\ \cite[VI.1.2]{rc-book}.
     Thus Proposition~\ref{cones.unique.pr} is the expected result for such cones.

     By contrast,  if $X$ is Fano and klt, there can be many
     proper, birational morphisms  $\pi:Y\to  \hat R(X, \Delta)$
     such that the exceptional set $E$ is a klt, Fano variety with negative normal bundle; see \cite {MR3187625, k-nqmmp}. 
     Then $E$ determines a filtration whose corresponding graded ring is a cone over $E$, so it looks like $\hat R(X, \Delta)$ could be isomorphic to a cone over $E$.

     Nonetheless, Proposition~\ref{cones.unique.pr} says that if $\auts^\circ(X, \Delta)$ is unipotent, then
     these other filtrations do not come from a grading.
     \end{rem}


   \def\cprime{$'$} \def\cprime{$'$} \def\cprime{$'$} \def\cprime{$'$}
  \def\cprime{$'$} \def\dbar{\leavevmode\hbox to 0pt{\hskip.2ex
  \accent"16\hss}d} \def\cprime{$'$} \def\cprime{$'$}
  \def\polhk#1{\setbox0=\hbox{#1}{\ooalign{\hidewidth
  \lower1.5ex\hbox{`}\hidewidth\crcr\unhbox0}}} \def\cprime{$'$}
  \def\cprime{$'$} \def\cprime{$'$} \def\cprime{$'$}
  \def\polhk#1{\setbox0=\hbox{#1}{\ooalign{\hidewidth
  \lower1.5ex\hbox{`}\hidewidth\crcr\unhbox0}}} \def\cdprime{$''$}
  \def\cprime{$'$} \def\cprime{$'$} \def\cprime{$'$} \def\cprime{$'$}
\providecommand{\bysame}{\leavevmode\hbox to3em{\hrulefill}\thinspace}
\providecommand{\MR}{\relax\ifhmode\unskip\space\fi MR }
\providecommand{\MRhref}[2]{%
  \href{http://www.ams.org/mathscinet-getitem?mr=#1}{#2}
}
\providecommand{\href}[2]{#2}

  \bigskip

  Princeton University, Princeton NJ 08544-1000, \

  \email{kollar@math.princeton.edu}

\end{document}